\documentclass{article}
\usepackage{bbm}
\usepackage{}
\usepackage{amssymb}
\usepackage{mathrsfs}
\usepackage{cite}

\usepackage{stmaryrd}
\usepackage{amsfonts}
\usepackage{amsmath}
\usepackage{bm}
\usepackage{indentfirst}
\usepackage{array}
\usepackage{arydshln}
\usepackage{graphicx}
\usepackage{listings}
\numberwithin{equation}{section}
\usepackage{geometry}
\usepackage{amsthm}
\title{The existence of full dimensional invariant tori for an almost-periodically forced nonlinear beam equation
\thanks{This work is supported by the NNSF of China (No.11971163)}}

\author {Shujuan Liu\thanks{ Corresponding author.}~~~~Guanghua Shi\\
\setcounter{footnote}{-1}
\noindent\footnote{~E-mail addresses: 846706840@smail.hunnu.edu.cn (S. Liu).}\\
{\small Key Laboratory of High Performance Computing and  Stochastic Information Processing,}\\
{\small Department of Mathematics, Hunan Normal University,}\\
{\small Changsha, Hunan 410083, China}
}
\date{}

\begin{document}\large
\maketitle
\textbf {Abstract.} In this paper, we prove the existence of full dimensional invariant tori for a non-autonomous, almost-periodically forced nonlinear beam equation with a periodic boundary condition via KAM theory.
\section{Introduction}
Presently, there have been many remarkable results on the beam equations via KAM theory, see\cite{J.Chang2015,JianshengGeng2003,L.H.Eliasson,JianshengGeng2006,JianshengGeng20061,X.Xu2009}.
 In these papers, the authors proved the existence of quasi-periodic solutions for  nonlinear beam equations, which is to say, the persistence of finite dimensional invariant tori for linear equations. In this paper, we will discuss the existence of almost-periodic solutions for a nonlinear beam equation. As to almost-periodic solutions, via KAM theory, Bourgain \cite{J.Bourgain2005} considered the nonlinear Schr\"{o}dinger equation $\sqrt{-1}u_{t}-u_{xx}+Mu+f(|u|^{2})u=0$ with a periodic boundary condition, where $M$ is a random Fourier multiplier. He proved that, for appropriate $M$, the above equation has invariant tori of full dimension. While for fixed $M$, by extracting parameters based on Birkhoff normal form, in \cite{JianshengGeng2012} it has been proved that the equation admits a family of small-amplitude full-dimensional invariant tori. Niu and Geng \cite{H.Niu2007} obtained almost periodic solutions for the case of higher dimensional beam equations.

\par The above equations do not depend on the forced terms. Physically, it means that there is no external force acting when the string is at rest, tending to distort its equilibrium of $u\equiv0$. As to the case with a forced term, Zhang-Si \cite{Zhang-Si2009} proved the existence of quasi-periodic solutions for the quasi-periodically forced nonlinear wave equation

\begin{equation}\label{z}
 u_{tt}=u_{xx}-\mu u-\varepsilon\phi(t)h(u),
\end{equation}
where $\mu>0$. If $\mu=0,$ the wave equation is completely resonant and we cannot directly extract parameters as usual. To overcome this difficulty, Yuan \cite{Yuan2005} has proved that (\ref{z}) with $\mu=0$ and $\phi(t)\equiv1$ has some special solutions depending only on the space variable $x$, and regarded such solutions as `parameters'. Along Yuan's idea, Si \cite{J.Si2012} studied the quasi-periodic solutions for non-autonomous quasi-periodically forced nonlinear wave equation
$$u_{tt}-u_{xx}+\phi(t,\varepsilon)u^3=0$$ with periodic boundary conditions via infinite-dimensional KAM theory developed by Kuksin \cite{Kusin.S.B.1993}. Besides, by the method of Lyapunov-Schmidt decomposition, Berti and Procesi \cite{M.Berti} showed the existence of small amplitude quasi-periodic solutions with two frequencies $\omega=(\omega_{1},\omega_{2})=(\omega_{1},1+\varepsilon)$ for completely resonant wave equation with periodic forcing
 $$\left\{
     \begin{array}{ll}
       u_{tt}-u_{xx}+f(\omega_{1}t,u)=0,\\
       u(t,x)=u(t,x+2\pi),\\
     \end{array}
   \right.$$ where the nonlinear forced term $f(\omega_{1}t,u)=a(\omega_{1}t)u^{2d+1}+O(u^{2d+2}),~d\in\mathbb{N}=\{1,2,\cdots\},$ is $2\pi/\omega_{1}$-periodic in time. By the similar method in \cite{Yuan2005}, Rui and Si \cite{J.Rui2014} turned the inhomogeneous  Schr\"{o}dinger equation $\sqrt{-1}u_{t}-u_{xx}+mu+\phi(t)|u|^2u=\varepsilon g(t)$ into a complex ordinary differential equation uniformly in space variables and a PDE with zero equilibrium point, then constructed the invariant tori or quasi-periodic solutions of the PDE by KAM method.
 For the existence of quasi-periodic solutions for non-autonomous PDEs, there are a couple of other references, see \cite{Y.Wang2012,Y.Shi2017,M.Zhang2016,L.Jiao2009}.
Up to now, there are only a few results about the existence of almost-periodic solutions for Hamiltonian partial differential equations (HPDEs) with almost-periodic forcing. Firstly, for the linear Schr\"{o}dinger equation with almost-periodic forcing
 $$\sqrt{-1}u_{t}-u_{xx}+mu+\psi(t)(a_{1}u+a_{2}\bar{u})=0,~~t\in\mathbb{R},~x\in\mathbb{T}^{1}$$under periodic boundary conditions, the existence of almost-periodic solutions is discussed by Rui, Liu and Zhang \cite{J.Rui2016}. Later on, based on reducibility via an improved KAM method, Rui-Liu \cite{J.Rui2017} focused on almost-periodic solutions for the linear wave equation with almost-periodic forcing
 $$u_{tt}=u_{xx}-mu-\psi(\omega t, x)u,~~t\in\mathbb{R},~~x\in[0,\pi]$$
 subject to periodic boundary conditions. Both of the above equations with almost-periodic forcing are linear. A natural question is that whether or not there are some almost-periodic solutions for the nonlinear HPDEs with almost-periodic forcing. The aim of this paper is to discuss the existence of almost-periodic solutions.
\par In this paper, we consider the following nonlinear beam equation with almost-periodic forcing
\begin{equation}
\begin{aligned}\label{eq1.1}
u_{tt}+u_{xxxx}+mu+\psi_{0}(\omega t)+\psi_{1}(\omega t)u+\psi_{2}(\omega t)u^2+\psi_{3}(\omega t)u^3=0,~~t\in\mathbb{R},~~x \in[0,2\pi],
\end{aligned}
\end{equation}
subject to the periodic boundary condition
\begin{equation}
\begin{aligned}\label{eq1.2}
u(t,0)=u(t,2\pi),~~t\in \mathbb{R},
\end{aligned}
\end{equation}
where $m>0,$ and $\psi_{l}(\omega t)~ (l=0,1,2,3)$ are almost-periodic in time.

\par Notice that (\ref{eq1.1}) is inhomogeneous if $\psi_{0}(\omega t)\not\equiv0$, therefore, it is easy to see that $u\equiv0$ is not a solution of (\ref{eq1.1})+(\ref{eq1.2}). Instead of the method in \cite{J.Rui2014} to deal with inhomogeneous terms, we directly begin with the initial system. Moreover, since the forced terms are almost-periodic, the main difficulty is to deal with the infinitely many frequencies. Following the idea of Rui an Liu \cite{J.Rui2017}, we add some conditions on the almost-periodic forcing such that we can choose properly finite frequencies at each KAM step to deal with small divisors.
\par To state our conditions on the almost-periodic forcing, we introduce some notations.

\par Consider the frequencies $\omega$ of the almost-periodic function $\psi_{l}(\omega t)~(l=0,1,2,3)$, and denote $\theta=\omega t\in\mathbb{T}^{\infty}$, where $\mathbb{T}^{\infty}$ is the infinite-dimensional torus. Let $\mathcal{O}$ denote the closed set $[0,1]^{\infty}$, that is $$\mathcal{O}=\{\xi:~\xi=(\xi_{1},\xi_{2},\cdots),~\xi_{j}\in[0,1],~j=1,2,\cdots\}. $$
$\psi_{l}(\omega t),~l=0,1,2,3$ are almost-periodic functions.
In this paper, we regard $\omega=(\omega_{1},\cdots)$ as parameters. Since the parameter set $\mathcal{O}$ is of infinite dimension, we explain the positive measure in the sense of Remark 1.2. Actually, it is enough to assume the parameter set is of finite dimension at every KAM iteration. We write
$$\omega=(\omega_{i_{1}},\omega_{i_{2}},\cdots,\omega_{i_{n}},\omega_{i_{n+1}},\cdots)=:(\omega^{n},\omega_{n}^{'})\in\mathcal{O}^{n}\times\mathcal{O}_{n}^{'}\equiv\mathcal{O},$$
where $(i_{1},i_{2},i_{3},\cdots)$ is a rearrangement of $(1,2,3,\cdots)$ and $\mathcal{O}^{n}=[0,1]^{n}$.
\par We choose a sequence $\{b_{v}\}$ satisfying $b_{0}=b\geq1,$ $b_{v+1}>b_{v},$ and $b_{v}\in\mathbb{Z}^{+}, v=0,1,\cdots$.
Let $$\mathcal{I}_{v}=\{i_{j}:j\leq b_{v},~j,~i_{j}\in\mathbb{Z}^{+}\},~v=0,1,2,\cdots,$$
thus, $\mathcal{I}_{\infty}:=\lim_{v\rightarrow\infty}\mathcal{I}_{v}=\mathbb{Z}^{+}.$\\
Denote
$$\omega_{1}^{b_{0}}=(\omega_{i_{1}},\cdots,\omega_{i_{b_{0}}}),~~~~\omega_{1}^{b_{v+1}}=(\omega_{i_{(b_{v}+1)}},\cdots,\omega_{i_{b_{(v+1)}}}),
~~\omega^{b_{v}}=(\omega_{i_{1}},\cdots,\omega_{i_{b_{v}}});$$
$$\theta_{1}^{b_{0}}=(\theta_{i_{1}},\cdots,\theta_{i_{b_{0}}}),~~~~\theta_{1}^{b_{v+1}}=(\theta_{i_{(b_{v}+1)}},\cdots,\theta_{i_{b_{(v+1)}}}),
~~\theta^{b_{v}}=(\theta_{i_{1}},\cdots,\theta_{i_{b_{v}}});$$
$$J_{1}^{b_{0}}=(J_{i_{1}},\cdots,J_{i_{b_{0}}}),~~~~J_{1}^{b_{v+1}}=(J_{i_{(b_{v}+1)}},\cdots,J_{i_{b_{(v+1)}}}),~~J^{b_{v}}=(J_{i_{1}},\cdots,J_{i_{b_{v}}}),$$
and let $[\psi]$ denotes the average of $\psi(\theta)$ on $\theta$, where $J_{i_{k}}$ is the action variable corresponding to angle variable $\theta_{i_{k}}$ later.
\par Throughout the paper, we assume that the small parameter $\varepsilon$ satisfies $0<\varepsilon\ll 1$ and that the following assumptions (\textbf{H1}) and (\textbf{H2}) hold.
\\
\\
  (\textbf{H1}) The functions $\psi_{l}(\omega t)$, $l=0,1,2,3$ are real analytic and almost-periodic in $t$ with frequencies $\omega$.\\
  \\
  (\textbf{H2}) For $0<\rho\leq1$, $\psi_{l}(\theta)=\sum_{j=0}^{\infty}\varepsilon^{(1+\rho)^{j}}\psi_{l}^{b_{j}}(\theta_{1}^{b_{j}}),$ $l=0,1,2,3$ are absolutely convergent, $\psi_{0}(\theta)\not\equiv0$ and there exists an absolute constant $C_{0}$ such that $$|\psi_{l}^{b_{j}}(\theta_{1}^{b_{j}})|\leq C_{0},~~|\partial_{\theta_{j^{'}}}\psi_{l}^{b_{0}}(\theta_{1}^{b_{0}})|\leq C_{0},~~j^{'}\in\mathcal{I}_{0},$$$$|\partial_{\theta_{j^{'}}}\psi_{l}^{b_{j}}( \theta_{1}^{b_{j}})|\leq C_{0},~~j^{'}\in\mathcal{I}_{j}\setminus\mathcal{I}_{j-1},~~l=0,1,2,3,~~j=1,2,\cdots,$$
where $|\cdot|$ denotes the sup-norm on $\mathbb{T}^{b_{j}}$ and $\partial_{\theta_{j^{'}}}f$ denotes the partial derivative of $f$ with respect to $\theta_{j^{'}}$.
\\
\textbf{Theorem 1.1}
Assume that the beam equation (\ref{eq1.1}) with periodic boundary condition (\ref{eq1.2}) satisfies the conditions (\textbf{H1}) and (\textbf{H2}). For $m>0$ and $0<\rho\leq1,$ there exists a positive measure Cantor-like subset $\mathcal{O}^{*}\subset\mathcal{O}$ such that for each $\omega=(\omega_{i_{1}},\omega_{i_{2}},\cdots)_{i_{j}\in\mathcal{I}_{\infty}}\in\mathcal{O}^{*},$ the beam equation (\ref{eq1.1})+(\ref{eq1.2}) has an almost-periodic solution of the form
$$u(t,x)=\sum_{j\geq0}\frac{q_{j}(\omega t)\cos(jx)}{\sqrt{\mu_{j}}},$$
where $\mu_{j}=\sqrt{j^4+m},$ $q_{j}(\omega t),~ j=0,1,\cdots$ are almost-periodic in $t$ with frequencies $\omega$ and $\mathop{\text{sup}}\limits_{\theta\in\mathbb{T}^{\infty}}\|q(\theta)\|_{\ell^{a,p+2}}=O(\varepsilon^{\frac{1}{2}-\frac{1}{8}\rho})$ with $p>0$. \\
\\
\textbf{Remark 1.2} Let the set $\mathcal{O}=[0,1]^{\infty}$ with probability measure. We say $\mathcal{O}^{*}$ is a set of large measure in $\mathcal{O}$ if there exists a real number $\gamma>0$ and $0<\varepsilon\ll1$ such that the following inequality holds
\begin{equation*}
  \text{meas}(\mathcal{O}\setminus\mathcal{O}^{*})\leq C\varepsilon^{\gamma},
\end{equation*}
where $\text{meas}$ is the standard probability measure on $[0,1]$ and $C>0$ is an absolute constant.\\\\
\textbf{Remark 1.3} Theorem 1.1 still holds for $\rho>1$. In this case, set $\varepsilon_{v}=\varepsilon^{\frac{1}{2}(1+\frac{1}{2\rho})^{v}}$, then we obtain a almost-periodic solution with the same form as that in Theorem 1.1 and $q$ satisfies $\mathop{\text{sup}}\limits_{\theta\in\mathbb{T}^{\infty}}\|q(\theta)\|_{\ell^{a,p+2}}=O(\varepsilon^{\frac{1}{2}-\frac{1}{8\rho}}).$
\\

Here is the outline of the proof of Theorem 1.1. Consider the perturbed beam equation
\begin{equation*}
u_{tt}+u_{xxxx}+mu+\psi_{0}(\omega t)+\psi_{1}(\omega t)u+\psi_{2}(\omega t)u^2+\psi_{3}(\omega t)u^3=0,~~t\in\mathbb{R},~~x \in[0,2\pi], \theta=\omega t.
\end{equation*}
It is easy to check that the above equation corresponds to a Hamiltonian equation
 \begin{equation}
\begin{aligned}\label{e}
\left\{
     \begin{array}{ll}
       \dot{u}=\frac{\partial H}{\partial v}=v\\
       \dot{v}=-\frac{\partial H}{\partial u}=-\big(u_{xxxx}+mu+\psi_{0}(\theta)+\psi_{1}(\theta)u+\psi_{2}(\theta)u^2+\psi_{3}(\theta)u^3\big)
     \end{array}
   \right.
\end{aligned}
\end{equation}
with the Hamiltonian
\begin{align*}
 H=&\frac{1}{2}\langle v,v\rangle+\frac{1}{2}\langle Au,u\rangle+\int_{0}^{2\pi}\psi_{0}(\theta)udx+\frac{1}{2}\int_{0}^{2\pi}\psi_{1}(\theta)u^2dx\\
&+\frac{1}{3}\int_{0}^{2\pi}\psi_{2}(\theta)u^3dx+\frac{1}{4}\int_{0}^{2\pi}\psi_{3}(\theta)u^4dx,
\end{align*}
where $A=\partial_{xxxx}+m$. By Fourier transformation, the above Hamiltonian can be turned into the infinite dimensional Hamiltonian system with Hamiltonian
\begin{equation*}
  H=N+P(\theta, z, \bar{z},\omega),
\end{equation*}
 with $z$, $\bar{z}\in\ell^{a,p}(\mathbb{C}):=\{z=(z_{0},z_{1},z_{2},\cdots):z_{i}\in\mathbb{C},|z_{0}|^2+\sum_{j\geq1}|z_{j}|^2j^{2p}e^{2ja}<\infty\}$ for $a>0$ and $p>0$.
 Since the forced term $\psi_{l}(\theta)(l=0,1,2,3)$ is an almost-periodic function with infinite frequencies $\omega=(\omega_{1},\cdots)$, we have to face the problem that how to treat infinite frequencies in the procedure of constructing almost-periodic solutions.
 But through observation, note that the forced terms still have some good `properties' which take part in proving Theorem 1.1.  $\psi_{l}(\theta)$ can be written as the sum of infinite functions with each one depending on a finite dimensional vector $\theta_{1}^{j} (j=b_{0},b_{1},\cdots)$ . In detail,
\begin{equation*}
  \psi_{l}(\theta)=\psi_{l}^{b_{0}}(\theta_{1}^{b_{0}})+\varepsilon^{(1+\rho)}\psi_{l}^{b_{1}}(\theta_{1}^{b_{1}})+\varepsilon^{(1+\rho)^{2}}\psi_{l}^{b_{2}}(\theta_{1}^{b_{2}})\cdots.
\end{equation*}
 According to the  property of $\psi_{l}(\theta)$ $(l=0,1,2,3)$, the perturbation $P(\theta,z,\bar{z})$ of  the Hamiltonian $H=N+P$ can also be written as
\begin{equation*}
  P(\theta,z,\bar{z})=\sum_{n\geq0}\varepsilon^{(1+\rho)^{n}}\tilde{P}(\theta_{1}^{b_{n}},z,\bar{z},\omega_{1}^{b_{n}}).
\end{equation*}
Therefore, we will construct the almost-periodic solution as follows. Firstly, we split the perturbation $P$ into two parts, that is
\begin{equation*}
  P(\theta,z,\bar{z})=\tilde{P}^{b_{0}}(\theta^{b_{0}},z,\bar{z},\omega^{b_{0}})+\sum_{n\geq1}\varepsilon^{(1+\rho)^{n}}\tilde{P}^{b_{n}}(\theta_{1}^{b_{n}},z,\bar{z},\omega_{1}^{b_{n}}),
\end{equation*}
where the former is the part depending only on $z,$ $\bar{z}$ and the finite dimensional vectors $\theta^{b_{0}}$, $\omega^{b_{0}}$, while, the latter is the part not depending on $\theta_{1}^{b_{0}}$. Secondly, spilt $\tilde{P}^{b_{0}}(\theta^{b_{0}},z,\bar{z},\omega^{b_{0}})$ into two parts,
\begin{equation*}
  \tilde{P}^{b_{0}}=(\tilde{P}^{b_{0}})^{low}+(\tilde{P}^{b_{0}})^{high},
  \end{equation*}
  where
  \begin{equation*}
    (\tilde{P}^{b_{0}})^{low}=\sum_{\gamma,\kappa\in\mathbb{N}^{\mathbb{N}},|\gamma|_{1}+|\kappa|_{1}\leq2}(\tilde{P}^{b_{0}})^{\gamma\kappa}(\theta^{b_{0}},\omega^{b_{0}})z^{\gamma}\bar{z}^{\kappa}
  \end{equation*}
is the part of low order, and
\begin{equation*}
    (\tilde{P}^{b_{0}})^{high}=\sum_{\gamma,\kappa\in\mathbb{N}^{\mathbb{N}},|\gamma|_{1}+|\kappa|_{1}\geq3}(\tilde{P}^{b_{0}})^{\gamma\kappa}(\theta^{b_{0}},\omega^{b_{0}})z^{\gamma}\bar{z}^{\kappa}
  \end{equation*}
  is the part of high order. Expand $(\tilde{P}^{b_{0}})^{\gamma\kappa}(\theta^{b_{0}},\omega^{b_{0}})$ into Fourier series
  \begin{equation*}
    (\tilde{P}^{b_{0}})^{\gamma\kappa}(\theta^{b_{0}},\omega^{b_{0}})=\sum_{k\in\mathbb{Z}^{b_{0}}}(\tilde{P}^{b_{0}})^{\gamma\kappa}_{k}e^{\sqrt{-1}\langle k,\theta^{b_{0}}\rangle}.
  \end{equation*}
   As the standard KAM procedure, we will remove all non-normalized terms
  \begin{equation*}
   \sum_{|\gamma|_{1}+|\kappa|_{1}\leq2;|k|+|\gamma-\kappa|\neq0}\sum_{k\in\mathbb{Z}^{b_{0}}}(\tilde{P}^{b_{0}})^{\gamma\kappa}_{k}e^{\sqrt{-1}\langle k,\theta^{b_{0}}\rangle}z^{\gamma}\bar{z}^{\kappa}
  \end{equation*}
in $(\tilde{P}^{b_{0}})^{low}$ by a symplectic transformation $\Phi_{0}$ . After the first step, we obtain the new Hamiltonian
\begin{equation*}
  H_{1}=N_{1}+P_{1}
\end{equation*}
where the new perturbation $P_{1}$ can be written in the form
\begin{equation*}
  P_{1}=\hat{P}_{1}(\theta^{b_{1}},z,\bar{z},\omega^{b_{1}})+\sum_{n\geq2}\varepsilon^{(1+\rho)^{n}}\tilde{P}^{b_{n}}(\theta_{1}^{b_{n}},z,\bar{z},\omega_{1}^{b_{n}})\circ\Phi_{0}.
\end{equation*}
It is easy to check that $\hat{P}_{1}(\theta^{b_{1}},z,\bar{z},\omega^{b_{1}})$ depends only on $z,$ $\bar{z}$ and the finite dimensional vectors $\theta^{b_{1}}$, $\omega^{b_{1}}$.
As the first step, we will split $\hat{P}_{1}$ into $\hat{P}_{1}^{low}$ and $\hat{P}_{1}^{high}$ and remove the non-normalized terms in $\hat{P}_{1}^{low}$. In the end, after infinite transformations, we obtain a non-degenerate normal form
\begin{equation*}
  H^{\infty}=\langle \omega,J\rangle+\langle \tilde{\Omega},z\bar{z}\rangle+\sum_{|\gamma|_{1}+|\kappa|_{1}\geq3}\hat{P}^{\gamma\kappa}(\theta,\omega)z^{\gamma}\bar{z}^{\kappa},
\end{equation*}
where $\tilde{\Omega}_{j}$ is close to the eigenvalue $\mu_{j}$ of operator $A$. Finally, basing on the above normal form, we get the almost-periodic solution in Theorem 1.1.\\
\indent In the process of removing the non-normalized terms, the following non-resonant conditions are needed,
\begin{equation*}\label{eq3.2}
  |\langle k,\omega^{b_{v}}\rangle+\langle l,\Omega_{v}\rangle|\geq \frac{\alpha_{v}\langle l\rangle_{2}}{(1+v^2)(|k|+1)^{2b_{v}+2}}
\end{equation*}
 for any $(k,l)\in \mathcal{Z}^{b_{v}}=\{(k,l)\neq0,|l|\leq2\}\subset \mathbb{Z}^{b_{v}}\times\ \mathbb{Z}^{\infty}$, where $v$ represents the iteration step.
 It is known that the standard non-resonant conditions depend on the dimension of torus. In the present paper, the dimension of torus increases with the iteration. Therefore, the non-resonant conditions depend on iteration steps.

\par The rest of this paper is organized as follows. In Section 2, we introduce the Hamiltonian setting of (\ref{eq1.1})+(\ref{eq1.2}). Sections 3--5 are devoted to the proof of Theorem 1.1. Some technical lemmas are listed in Appendix.
\\

\section{Hamiltonian setting}
The system (\ref{eq1.1})+(\ref{eq1.2}) can be written as a Hamiltonian system
\begin{equation}
\begin{aligned}\label{e}
\left\{
     \begin{array}{ll}
       \dot{u}=\frac{\partial H}{\partial v}=v\\
       \dot{v}=-\frac{\partial H}{\partial u}=-\big(u_{xxxx}+mu+\psi_{0}(\omega t)+\psi_{1}(\omega t)u+\psi_{2}(\omega t)u^2+\psi_{3}(\omega t)u^3\big)
     \end{array}
   \right.
\end{aligned}
\end{equation}
with the Hamiltonian
\begin{align*}
 H=&\frac{1}{2}\langle v,v\rangle+\frac{1}{2}\langle Au,u\rangle+\int_{0}^{2\pi}\psi_{0}(\omega t)udx+\frac{1}{2}\int_{0}^{2\pi}\psi_{1}(\omega t)u^2dx\\
&+\frac{1}{3}\int_{0}^{2\pi}\psi_{2}(\omega t)u^3dx+\frac{1}{4}\int_{0}^{2\pi}\psi_{3}(\omega t)u^4dx,
\end{align*}
where $A=\partial_{xxxx}+m$.
The eigenvalues of the operator $A$ with the periodic boundary condition are $\mu_{j}^2=j^4+m$, $j\in\mathbb{Z}$, corresponding eigenfunction $\phi_{j}(x)\in L^{2}[0,2\pi]$

$$\phi_{j}(x)=\left\{
     \begin{array}{ll}
       \frac{1}{\sqrt{\pi}}\cos jx,~~~~j>0,\\
       -\frac{1}{\sqrt{\pi}}\sin jx,~~j<0,\\
       \frac{1}{\sqrt{2\pi}},~~~~~~~~~~~j=0.\\
     \end{array}
   \right.$$
In order to avoid the double eigenvalues, we restrict ourselves to find some solutions which are even in $x$. $\{\phi_{j}:j\geq0\}$ is a complete orthogonal basis of a subspace in $L^2[0,2\pi]$.
\par We introduce coordinates $q=(q_{0},q_{1},q_{2},\cdots)$ and $\chi=(\chi_{0},\chi_{1},\chi_{2},\cdots)$ by the following relations
$$u(t,x)=\sum_{j\geq0}\frac{q_{j}(t)}{\sqrt{\mu_{j}}}\phi_{j}(x),~~~~v(t,x)=\sum_{j\geq0}\sqrt{\mu_{j}}\chi_{j}(t)\phi_{j}(x).$$
The coordinates are taken from some real Hilbert space
$$\ell^{a,p}=\ell^{a,p}(\mathbb{R}):=\{q=(q_{0},q_{1},q_{2},\cdots):q_{j}\in\mathbb{R},j\geq0\}$$ with norm $$\|q\|_{a,p}^{2}=|q_{0}|^2+\sum_{j\geq1}|q_{j}|^2j^{2p}e^{2ja}<\infty.$$
In the sequel, we assume that $a\geq0$ and $p>0$.
\par We introduce a pair of action-angle variables $(J,\theta)\in\mathbb{R}^{\infty}\times\mathbb{T}^{\infty}$ such that the Hamiltonian is autonomous. We then obtain the Hamiltonian
\begin{align*}
 H=&\langle\omega,J\rangle+\frac{1}{2}\sum_{j\geq0}\mu_{j}(\chi_{j}^2+q_{j}^2)\\
&+\int_{0}^{2\pi}\psi_{0}(\theta)\sum_{j\geq0}\frac{q_{j}(t)}{\sqrt{\mu_{j}}}\phi_{j}(x)dx
+\frac{1}{2}\int_{0}^{2\pi}\psi_{1}(\theta)\Big(\sum_{j\geq0}\frac{q_{j}(t)}{\sqrt{\mu_{j}}}\phi_{j}(x)\Big)^2dx\\
&+\frac{1}{3}\int_{0}^{2\pi}\psi_{2}(\theta)\Big(\sum_{j\geq0}\frac{q_{j}(t)}{\sqrt{\mu_{j}}}\phi_{j}(x)\Big)^3dx+\frac{1}{4}\int_{0}^{2\pi}\psi_{3}(\theta)\Big(\sum_{j\geq0}\frac{q_{j}(t)}{\sqrt{\mu_{j}}}\phi_{j}(x)\Big)^4dx\\
\end{align*}
with equations of motions
$$\dot{\theta}=\omega,~~\dot{J}=-\frac{\partial H}{\partial\theta},~~\dot{q}_{j}=\frac{\partial H}{\partial \chi_{j}},~~\dot{\chi}_{j}=-\frac{\partial H}{\partial q_{j}},~~j\geq0$$
 with respect to the symplectic structure $d\theta\wedge dJ+\sum_{j\geq0}dq_{j}\wedge d\chi_{j}.$  We introduce complex coordinates
$$z_{j}=\frac{1}{\sqrt{2}}(q_{j}-\sqrt{-1}\chi_{j}),~~~~~\bar{z}_{j}=\frac{1}{\sqrt{2}}(q_{j}+\sqrt{-1}\chi_{j}),~~~~j\geq0,$$
which are in the complex Hilbert space $\ell^{a,p}(\mathbb{C})$. Here
$$\ell^{a,p}(\mathbb{C}):=\{z=(z_{0},z_{1},z_{2},\cdots):z_{j}\in\mathbb{C},j\geq0\}$$ with  finite norm $$\|z\|_{a,p}^{2}=|z_{0}|^2+\sum_{j\geq1}|z_{j}|^2j^{2p}e^{2ja}<\infty.$$
We then obtain the Hamiltonian
\begin{equation}
\begin{aligned}\label{eq2.1}
H=\sum_{j=1}^{\infty}\omega_{i_{j}}J_{i_{j}}+\sum_{j\geq0}\mu_{j}z_{j}\bar{z}_{j}+G_{1}+G_{2}+G_{3}+G_{4}
\end{aligned}
\end{equation}
and the symplectic structure $d\theta\wedge dJ+\sqrt{-1}\sum_{j\geq0}d\bar{z}_{j}\wedge dz_{j},$ where

\begin{eqnarray}
  &&G_{1}=\sqrt{\pi}\psi_{0}(\theta)\frac{z_{0}+\bar{z}_{0}}{\sqrt{\mu_{0}}},~~~~G_{2}=\frac{1}{4}\psi_{1}(\theta)\sum_{j\geq0}\frac{(z_{j}+\bar{z}_{j})^2}{\mu_{j}},\label{G12}\\
  &&G_{3}=\frac{\sqrt{2}}{12}\psi_{2}(\theta)\sum_{i\pm j\pm l=0,i,j,l\geq0}\frac{G_{ijl}^3}{\sqrt{\mu_{i}\mu_{j}\mu_{l}}}(z_{i}+\bar{z}_{i})(z_{j}+\bar{z}_{j})(z_{l}+\bar{z}_{l}),\label{G3}  \\
  &&G_{4}=\frac{1}{16}\psi_{3}(\theta)\sum_{i\pm j\pm k \pm l=0,i,j,k,l\geq0}\frac{G_{ijkl}^4}{\sqrt{\mu_{i}\mu_{j}\mu_{k}\mu_{l}}}(z_{i}+\bar{z}_{i})(z_{j}+\bar{z}_{j})(z_{k}+\bar{z}_{k})(z_{l}+\bar{z}_{l}),\label{G4}
\end{eqnarray}
 $$G_{ijl}^3=\int_{0}^{2\pi}\phi_{i}\phi_{j}\phi_{l}dx,~~~~G_{ijkl}^4=\int_{0}^{2\pi}\phi_{i}\phi_{j}\phi_{k}\phi_{l}dx.$$
It is easy to prove that $G_{ijl}^3=0$ unless $i\pm j\pm l=0$, and $G_{ijkl}^4=0$ unless $i\pm j\pm k\pm l=0.$

\par According to the assumption (\textbf{H2}) on $\psi_{l}({\theta})~(l=0,1,2,3)$, we  split $G$ with $G=G_{1}+G_{2}+G_{3}+G_{4}$ into infinitely many parts, that is $G=\sum_{n\geq0}\varepsilon^{(1+\rho)^{n}}\tilde{P}^{b_{n}}$, where $\tilde{P}^{b_{n}}$ only depends on $z$, $\bar{z}$, $\omega^{b_{n}}$ and angle variables $\theta_{1}^{b_{n}}$. Then we just need to treat finitely many frequencies at each step. More precisely, we write (\ref{eq2.1}) as
\begin{equation}\label{eqk}
  H=\langle \omega,J\rangle+\sum_{j\geq0}\mu_{j}z_{j}\bar{z}_{j}+\sum_{n\geq0}\varepsilon^{(1+\rho)^{n}}\tilde{P}^{b_{n}}
\end{equation}
with
\begin{equation}
\begin{aligned}\label{Pn}
\tilde{P}^{b_{n}}=&\frac{\sqrt{\pi}}{m^{\frac{1}{4}}}\psi_{0}^{b_{n}}(\theta_{1}^{b_{n}})(z_{0}+\bar{z}_{0})+\frac{1}{4}\psi_{1}^{b_{n}}(\theta_{1}^{b_{n}})\sum_{j\geq0}\frac{1}{\mu_{j}}(z_{j}^{2}+2z_{j}\bar{z}_{j}+\bar{z}_{j}^{2})\\
&+\frac{\sqrt{2}}{12}\psi_{2}^{b_{n}}(\theta_{1}^{b_{n}})\sum_{i\pm j\pm l=0,i,j,l\geq0}\frac{G_{ijl}^3}{\sqrt{\mu_{i}\mu_{j}\mu_{l}}}(z_{i}+\bar{z}_{i})(z_{j}+\bar{z}_{j})(z_{l}+\bar{z}_{l})\\
&+\frac{1}{16}\psi_{3}^{b_{n}}(\theta_{1}^{b_{n}})\sum_{i\pm j\pm k \pm l=0,i,j,k,l\geq0}\frac{G_{ijkl}^4}{\sqrt{\mu_{i}\mu_{j}\mu_{k}\mu_{l}}}(z_{i}+\bar{z}_{i})(z_{j}+\bar{z}_{j})(z_{k}+\bar{z}_{k})(z_{l}+\bar{z}_{l}),
 \end{aligned}
 \end{equation}
where $\psi_{l}^{b_{n}}(\theta_{1}^{b_{n}}),~l=0,1,2,3$ are defined in the assumption (\textbf{H2}).\\
\par  We define the Hamiltonian vector field of a Hamiltonian $Q$, $$X_{Q}=(Q_{J},-Q_{\theta},-\sqrt{-1}Q_{z},\sqrt{-1}Q_{\bar{z}}).$$
To obtain the analyticity of $X_{G}$, it is convenient to introduce coordinates $$w=(\cdots,w_{-1},w_{-0},w_{0},w_{1},\cdots)\in\ell_{b}^{a,p}$$ by setting $z_{j}=w_{j}$, $\bar{z}_{j}=w_{-j}$, where $\ell_{b}^{a,p}$ consists of all bi-infinite sequences with finite norm $$\|w\|_{a,p}^{2}=|w_{0}|^2+|w_{-0}|^2+\sum_{|j|\geq1}^{\infty}|w_{j}|^2|j|^{2p}e^{2|j|a}.$$
Substituting $z_{j}=w_{j}$ and $\bar{z}_{j}=w_{-j}$ into (\ref{G12}), (\ref{G3}) and (\ref{G4}), then we obtain
\begin{eqnarray*}
G_{1}&=&\sqrt{\pi}\psi_{0}(\theta)\frac{w_{0}+w_{-0}}{\sqrt{\mu_{0}}},~~G_{2}=\frac{1}{4}\psi_{1}(\theta)\sum_{j\geq0}\frac{(w_{j}+w_{-j})^2}{\mu_{j}}, \\
G_{3}&=&\frac{\sqrt{2}}{12}\psi_{2}(\theta)\sum_{i\pm j\pm l=0,i,j,l\in\mathbb{Z}}\frac{G_{|i||j||l|}^3}{\sqrt{\mu_{i}\mu_{j}\mu_{l}}}w_{i}w_{j}w_{l},  \\
G_{4}&=&\frac{1}{16}\psi_{3}(\theta)\sum_{i\pm j\pm k \pm l=0,i,j,k,l\in\mathbb{Z}}\frac{G_{|i||j||k||l|}^4}{\sqrt{\mu_{i}\mu_{j}\mu_{k}\mu_{l}}}w_{i}w_{j}w_{k}w_{l}.
\end{eqnarray*}

\indent \textbf {Lemma 2.1}([22,Lemma 2])
For $a\geq0$ and $p>\frac{1}{2},$ the space $\ell^{a,p}_{b}$ is a Hilbert algebra with respect to convolution of sequences, and
$$\|q*h\|_{a,p}\leq c\|q\|_{a,p}\|h\|_{a,p}$$
with a constant $c$ depending only on $p$.\\
\\
\indent \textbf {Lemma 2.2}
For $a\geq0,$ and $p>0,$ the Hamiltonian vector field $X_{G}$ is real analytic as a map from some neighbourhood of the origin in $\ell_{b}^{a,p}$ into $\ell_{b}^{a,p+2},$ with $\|X_{G}\|_{a,p+2}\leq C\varepsilon$ uniformly in $\theta\in\mathbb{T}^{\infty}$, where $C$ is a positive constant.\\

\textbf{Proof} Set $\tilde{w}_{j}=\frac{1}{\sqrt{\mu_{j}}}(|w_{j}|+|w_{-j}|).$ By the assumption (\textbf{H2}) and $G_{ijkl}^4=\int_{0}^{2\pi}\phi_{i}\phi_{j}\phi_{k}\phi_{l}dx$, we have
\begin{equation*}
\begin{aligned}
\big|\partial_{w_{l}}G_{4}\big|
&\leq\frac{1}{16}|\psi_{3}(\theta)|\sum_{i\pm j\pm k =l,i,j,k,l\in\mathbb{Z}}\frac{\big|G_{|i||j||k||l|}^4\big|}{\sqrt{\mu_{i}\mu_{j}\mu_{k}\mu_{l}}}|w_{i}w_{j}w_{k}|\\
&\leq\frac{C\varepsilon}{\sqrt{\mu_{l}}}\sum_{i+j+k=l,i,j,k,l\in\mathbb{Z}}\tilde{w}_{i}\tilde{w}_{j}\tilde{w}_{k}\\
&=\frac{C\varepsilon}{\sqrt{\mu_{l}}}(\tilde{w}\ast \tilde{w}\ast \tilde{w})_{l}.
\end{aligned}
\end{equation*}
If $w\in\ell_{b}^{a,p}$, then $\tilde{w}\in\ell_{b}^{a,p+1}$, and for $p>0,$ the latter is a Hilbert algebra by Lemma 2.1. Therefore, $\tilde{w}\ast \tilde{w}\ast \tilde{w}$ also belongs to $\ell_{b}^{a,p+1},$ and $(G_{4})_{w}\in\ell_{b}^{a,p+2}$ with
\begin{equation*}
\begin{aligned}
\|(G_{4})_{w}\|_{a,p+2}\leq C\varepsilon\|\tilde{w}\ast \tilde{w}\ast \tilde{w}\|_{a,p+1}\leq C\varepsilon\|w\|_{a,p}^3.
\end{aligned}
\end{equation*}
Similarily,
$$\|(G_{3})_{w}\|_{a,p+2}\leq C\varepsilon\|w\|_{a,p}^2.$$
We also obtain that
\begin{equation*}\label{}
\begin{aligned}
\|(G_{1})_{w}\|_{a,p+2}^{2}=\big|{\partial_{w_{0}}G_{1}}\big|^{2}+\big|\partial_ {w_{-0}}G_{1}\big|^{2}=2\Big|\frac{\sqrt{\pi}\psi_{0}(\theta)}{\sqrt{\mu_{0}}}\Big|^{2}\leq C\varepsilon^{2},
\end{aligned}
\end{equation*}
and
\begin{equation*}\label{}
\begin{aligned}
\big|\partial_{w_{j}}G_{2}\big|=\frac{1}{2}|\psi_{1}(\theta)|\frac{|w_{j}+w_{-j}|}{\mu_{j}}\leq C\varepsilon\frac{\tilde{w}_{j}}{\sqrt{\mu_{j}}}.
\end{aligned}
\end{equation*}
Hence,
\begin{equation*}
\begin{aligned}
\|(G_{2})_{w}\|_{a,p+2}\leq C\varepsilon\|\tilde{w}\|_{a,p+1}\leq C\varepsilon\|w\|_{a,p}.
\end{aligned}
\end{equation*}
The proof of Lemma 2.2 is complete.   ~~~~~$\Box$

Let
\begin{equation}
\begin{aligned}\label{eqG}
H=N+P=\sum_{j=1}^{\infty}\omega_{i_{j}}J_{i_{j}}+\sum_{j\geq0}\Omega_{j}(\omega)z_{j}\bar{z}_{j}+P(\theta,z,\bar{z},\omega)
\end{aligned}
\end{equation}
be a Hamiltonian defined on a phase space $\mathcal{P}^{a,p}:=\mathbb{T}^{\infty}\times\mathbb{R}^{\infty}\times\ell^{a,p}\times\ell^{a,p}$, with the normal form $N=\sum_{j=1}^{\infty}\omega_{i_{j}}J_{i_{j}}+\sum_{j\geq0}\Omega_{j}(\omega)z_{j}\bar{z}_{j},$ and the perturbation $P(\theta,z,\bar{z},\omega)=\sum_{n\geq0}\varepsilon^{(1+\rho)^{n}}\tilde{P}^{b_{n}}$,
where $\omega=(\omega_{i_{1}},\omega_{i_{2}},\cdots)$ and $\Omega(\omega)=(\Omega_{0},\Omega_{1},\cdots)$ represent respectively, the tangent and normal frequencies.
Assume $P$ is analytic with respect to $\theta,~z,~\bar{z}$ and Lipschitz continuous in $\omega\in\mathcal{O}$. Note that, in this paper, we regard the tangent frequencies $\omega$ as parameters. When the perturbation vanishes, it is clear that $\mathcal{T}_{0}^{\infty}:=\mathbb{T}^{\infty}\times\{0\}\times\{0\}\times\{0\}$ depending on $\omega$ are infinitely dimensional invariant tori. Whether can these tori persist if $P$ is sufficiently small? In the following, we will prove that most of them can survive by the KAM method.
Firstly, we introduce some norms and notations.


The complex neighbourhood of torus $\mathcal{T}_{0}^{\infty}$ is defined by
$$D(s,r): |\text{Im}\theta|< s,~~|J|< r^{2},~~\|z\|_{a,p}<r,~~\|\bar{z}\|_{a,p}< r,$$
where $|\cdot|$ denotes the sup-norm for complex vectors, and weighted phase space norms are defined by

\begin{equation}\label{eq3.35}
|W|_{r}=|W|_{r,a,p+2}=|X|+\frac{1}{r^2}|Y|+\frac{1}{r}\|U\|_{a,p+2}+\frac{1}{r}\|V\|_{a,p+2},
\end{equation}
for $W=(X,Y,U,V)\in\mathcal{P}^{a,p+2}$.

Furthermore, we assume that the Hamiltonian vector field $X_{P}$ is real analytic on $D(s,r)$ for some positive $s$, $r$ uniformly in $\omega\in\mathcal{O}$ with finite norm $|X_{P}|_{r,D(s,r)\times\mathcal{O}}=\sup\limits_{D(s,r)\times\mathcal{O}}|X_{P}|_{r}$, and that the same holds for its Lipschitz semi-norm
$$|X_{P}|^{lip}_{r,\mathcal{O}}=\sup\limits_{\omega,\tilde{\omega}\in\mathcal{O},\omega\neq\tilde{\omega}}\frac{|\Delta_{\omega\tilde{\omega}}X_{P}|_{r}}{|\omega-\tilde{\omega}|},
~~~~|X_{P}|^{lip}_{r,D(s,r)\times\mathcal{O}}=\sup\limits_{D(s,r)}|X_{P}|^{lip}_{r,\mathcal{O}},$$
where $\Delta_{\omega\tilde{\omega}}X_{P}=X_{P}(\cdot,\omega)-X_{P}(\cdot,\tilde{\omega})$. Fixing $-2\leq\delta\leq0$, the Lipschitz semi-norm of the frequencies $\Omega(\omega)$ are defined by
$$
|\Omega|^{lip}_{-\delta,\mathcal{O}}=\sup\limits_{\omega,\tilde{\omega}\in\mathcal{O},\omega\neq\tilde{\omega}}\sup\limits_{j\geq0}\frac{j^{-\delta}|\Delta_{\omega\tilde{\omega}}\Omega_{j}|}{|\omega-\tilde{\omega}|}.$$
\par For $\lambda\geq0,$ define
$$|X_{P}|_{r,*}^{\lambda}=|X_{P}|_{r,*}+\lambda|X_{P}|_{r,*}^{lip},$$
where $*$ represents a set of variables (for example, $D(s,r)\times\mathcal{O}$), the symbol `$\lambda$' will always be used in this role and never have the meaning of exponentiation. Moreover, we introduce the notations   $$\langle l\rangle_{2}=\max(1,\big|\sum j^{2}l_{j}\big|),~~~~\mathcal{Z}^{b_{v}}=\{(k,l)\neq0,|l|\leq2\}\subset \mathbb{Z}^{b_{v}}\times\ \mathbb{Z}^{\infty}.$$

\section{Iteration lemma and its proof}
To state and prove the iterative lemma, we introduce some iterative constants and notations. Let $\varepsilon$, $s$, $r$ and $\rho$ be positive. Let $ v\geq0$ be the $v$-th KAM step, and set
\begin{itemize}
  \item [1.] $\varepsilon_{0}=\varepsilon^{\frac{1}{2}}$ and $\varepsilon_{v+1}=\varepsilon_{v}^{1+\frac{1}{2}\rho},$ $0<\rho\leq1;$
  \item [2.] $\alpha_{v}=\varepsilon_{v}^{\frac{1}{18}\rho},$ $M_{v}+1=(M_{1}+1)(2-2^{-v+1})~(v\geq1)$, $M_{0}=0$, $M_{1}=\varepsilon_{0}^{1-\frac{1}{4}\rho},$ $\lambda_{v}=\frac{\alpha_{v}}{M_{v}+1};$
  \item [3.]$\sigma_{v+1}=\frac{\sigma_{v}}{2},$ $s_{v+1}=s_{v}-6\sigma_{v},$ $s_{0}=s$ as initial value, fix $\sigma_{0}=s_{0}/24\leq1/20$ so that $s_{0}>s_{1}>\cdots\geq s_{0}/2;$
  \item [4.] $r_{v}=(1-\tau_{v})r_{0},$ with $\tau_{0}=0,$ $r_{0}=r$, $\tau_{v}=(1^{-2}+\cdot\cdot\cdot +v^{-2})/(2\sum^{\infty}_{j=1}j^{-2})$  $(v\geq1)$,
and $d_{v}=\frac{1}{4}(r_{v}-r_{v+1})=r_{0}/[8(v+1)^2\sum_{j=1}^{\infty}j^{-2}]$, $r_{0}>r_{1}>\cdots\geq r_{0}/2;$
  \item [5.]\begin{equation*}
\begin{aligned}
D_{v}=D(s_{v},r_{v})=\{&(\theta,J,z,\bar{z})\in\mathbb{C}^{\infty}/2\pi\mathbb{Z}^{\infty}\times\mathbb{C}^{\infty}\times \ell^{a,p}\times\ell^{a,p}:\\
&|\text{Im}\theta|< s_{v},~~|J|< r_{v}^{2},~~\|z\|_{a,p}<r_{v},~~\|\bar{z}\|_{a,p}< r_{v}\};
\end{aligned}
\end{equation*}

\item [6.]

$$\mathcal{O}_{v}^{\ast}=\{\omega:~\omega=(\omega_{i_{1}},\cdots,\omega_{i_{b_{v}}},\omega_{i_{(b_{v}+1)}},\dots)=:(\omega^{b_{v}},\omega_{b_{v}}^{'})\in\mathcal{O}_{*}^{b_{v}}\times\mathcal{O}_{b_{v}}^{'},~~i_{j}\in\mathbb{Z}^{+}\},$$
where $\mathcal{O}_{b_{v}}^{'}$ is the closed set of sequences $\omega_{b_{v}}^{'}=(\omega_{i_{(b_{v}+1)}},\omega_{i_{(b_{v}+2)}},\dots)$ with $\omega_{j}\in[0,1],$ $j=i_{(b_{v}+1)},\cdots,$
and
 $$\mathcal{O}_{*}^{b_{v}}=\mathcal{O}^{v}\backslash\big(\mathop{\bigcup}\limits_{k,l}\mathcal{R}_{kl}^{v}\big),$$
 $$\mathcal{O}^{v}=\{\omega^{b_{v}}:(\omega^{b_{v}},\omega_{b_{v}}^{'})\in\mathcal{O}_{v-1}^{*}\}\subset[0,1]^{b_{v}},~~~~~~~~\mathcal{O}_{-1}^{*}=\mathcal{O},$$
$$\mathcal{R}^{v}_{kl}=\big\{\omega^{b_{v}}\in\mathcal{O}^{v}:|\langle k,\omega^{b_{v}}\rangle+\langle l,\Omega_{v}\rangle|< \frac{\alpha_{v}\langle l\rangle_{2}}{\big(1+v^2\big)(|k|+1)^{2b_{v}+2}}\big\},~~~~(k,l)\in\mathcal{Z}^{b_{v}}.$$

\end{itemize}

\subsection{Iterative lemma}

We have obtained the Hamiltonian (\ref{eqk}) of (\ref{e}), which is of the form (\ref{eqG}) with the normal frequencies $\Omega=(\mu_{0},\mu_{1},\mu_{2},\cdots)$ and the perturbation $P=\sum_{n\geq0}\varepsilon^{(1+\rho)^n}\tilde{P}^{b_{n}}$. By Lemma 2.2 and the assumption (\textbf{H1}), it follows that $X_{G}$ is real analytic in $D(s,r)$ for some positive $s,$ $r$ uniformly in $\omega\in\mathcal{O}$, Lipschitz continuous in $\omega\in\mathcal{O}$ and $|\Omega|_{-\delta,\mathcal{O}}^{lip}=0.$
\par Assume that at the $v$-th step of scheme, a Hamiltonian
\begin{equation*}
  H_{v}=N_{v}+P_{v}=N_{v}+\hat{P}^{v}(\theta^{b_{v}},z,\bar{z},\omega^{b_{v}})+\sum_{n\geq v+1}\varepsilon^{(1+\rho)^{n}}\tilde{P}^{b_{n}}(\theta_{1}^{b_{n}},z,\bar{z},\omega_{1}^{b_{n}})\circ \Phi_{0} \circ \Phi_{1} \circ \cdots \circ \Phi_{v-1},
\end{equation*}
is considered as a small perturbation of some normal form $N_{v}.$ Split $\hat{P}_{v}$ into two parts, that is $\hat{P}_{v}=(\hat{P}_{v})^{low}+(\hat{P}_{v})^{high},$ where the low-degree terms $(\hat{P}_{v})^{low}$, denoted by $\varepsilon_{v}R^{b_{v}}$ in the following, is defined by
\begin{equation*}\label{Rv}
\varepsilon_{v}R^{b_{v}}=(\hat{P}_{v})^{low}=\sum_{|\gamma|_{1}+|\kappa|_{1}\leq2}\hat{P}_{v}^{\gamma\kappa}(\theta^{b_{v}},\omega^{b_{v}})z^{\gamma}\bar{z}^{\kappa},
\end{equation*}
\begin{equation*}
\begin{aligned}
R^{b_{v}}=&R^{00b_{v}}+\langle R^{10b_{v}},z\rangle+\langle R^{01b_{v}},\bar{z}\rangle+\langle R^{20b_{v}}z,z\rangle+\langle R^{11b_{v}}z,\bar{z}\rangle+\langle R^{02b_{v}}\bar{z},\bar{z}\rangle\\
=&\sum_{k\in\mathbb{Z}^{b_{v}}}R^{00b_{v}}_{k}e^{\sqrt{-1}\langle k,\theta^{b_{v}}\rangle}+\sum_{j\geq0, k\in\mathbb{Z}^{b_{v}}}(R_{kj}^{10b_{v}}z_{j}+R_{kj}^{01b_{v}}\bar{z}_{j})e^{\sqrt{-1}\langle k,\theta^{b_{v}}\rangle}\\
&+\sum_{i,j\geq0, k\in\mathbb{Z}^{b_{v}}}(R_{kij}^{20b_{v}}z_{i}z_{j}+R_{kij}^{11b_{v}}z_{i}\bar{z}_{j}+R_{kij}^{02b_{v}}\bar{z}_{i}\bar{z}_{j})e^{\sqrt{-1}\langle k,\theta^{b_{v}}\rangle}
\end{aligned}
\end{equation*}
and the high-degree terms $(\hat{P}_{v})^{high}$ of $\hat{P}_{v}$ by
$$(\hat{P}_{v})^{high}=\sum_{|\gamma|_{1}+|\kappa|_{1}\geq3}\hat{P}_{v}^{\gamma\kappa}(\theta^{b_{v}},\omega^{b_{v}})z^{\gamma}\bar{z}^{\kappa}.$$
The product $z^{\gamma}\bar{z}^{\kappa}$ denotes $\prod_{n}z_{n}^{\gamma_{n}}\bar{z}_{n}^{\kappa_{n}}$, where $\gamma=(\gamma_{0},\gamma_{1}\cdots,\gamma_{n},\cdots),$ $\kappa=(\kappa_{0},\kappa_{1},\cdots,\kappa_{n},\cdots)$ with finitely many non-zero components and $\gamma_{n},$ $\kappa_{n}\in\mathbb{N}$.\\

\textbf{Lemma 3.1}
Suppose that $H_{v}=N_{v}+P_{v}~(v\geq0)$ is given on $D_{v}\times\mathcal{O}_{v}^{*}$, where$$N_{v}=\langle\omega^{b_{v}},J^{b_{v}}\rangle+\sum_{n\geq v+1}\langle\omega_{1}^{b_{n}},J_{1}^{b_{n}}\rangle+\sum_{j\geq0}\Omega_{vj}(\omega^{b_{v-1}})z_{j}\bar{z}_{j}$$
is a normal form satisfying
\begin{equation*}\label{eq3.2}
  |\langle k,\omega^{b_{v}}\rangle+\langle l,\Omega_{v}\rangle|\geq \frac{\alpha_{v}\langle l\rangle_{2}}{(1+v^2)(|k|+1)^{2b_{v}+2}}, ~~~~(k,l)\in\mathcal{Z}^{b_{v}},
\end{equation*}

\begin{equation}\label{eq3.1}
|\Omega_{v}|^{lip}_{-\delta,\mathcal{O}_{v}^{\ast}}\leq M_{v},
\end{equation}
 \begin{equation}\label{eq3.4}
\Omega_{0j}=\mu_{j}=\sqrt{j^4+m},~~\Omega_{vj}=\Omega_{0j}+\sum_{\tilde{s}=0}^{v-1}\varepsilon_{\tilde{s}}[B_{jj}^{11b_{\tilde{s}}}],~~\Omega_{v}=(\Omega_{v0},\Omega_{v1},\Omega_{v2},\cdots), ~~~~v\geq1
\end{equation}
with $B_{jj}^{11b_{v}}$ defined in Section 3.2 and $P_{v}$ satisfies
\begin{equation*}
\begin{aligned}\label{eq3.21}
P_{v}=\hat{P}_{v}+\sum_{n\geq v+1}\varepsilon^{(1+\rho)^{n}}\tilde{P}^{b_{n}}(\theta_{1}^{b_{n}},z,\bar{z},\omega_{1}^{b_{n}})\circ \Phi_{0} \circ \Phi_{1} \circ \cdots \circ \Phi_{v-1}
\end{aligned}
\end{equation*}
with
\begin{equation*}\label{}
  \hat{P}_{v}=\varepsilon_{v}R^{b_{v}}(\theta^{b_{v}},z,\bar{z},\omega^{b_{v}})+\hat{P}_{v}^{high}(\theta^{b_{v}},z,\bar{z},\omega^{b_{v}}),
\end{equation*}

\begin{equation}\label{Ph}
|X_{\hat{P}_{v}^{high}}|_{r_{v},D(s_{v},r_{v})\times\mathcal{O}_{v}^{*}}^{\lambda_{v}}\leq\varepsilon_{0}+\sum_{j=1}^{v}\varepsilon_{j}^{\frac{1}{2}}
\end{equation}
and
\begin{equation}\label{eqi}
 |X_{R^{b_{v}}}|_{r_{v},D(s_{v},r_{v})\times\mathcal{O}_{v}^{*}}^{\lambda_{v}}\leq\frac{1}{2},~~~|X_{P_{v}-\hat{P}_{v}^{high}}|_{r_{v},D(s_{v},r_{v})\times\mathcal{O}_{v}^{*}}^{\lambda_{v}}\leq \varepsilon_{v}.
\end{equation}

Then there exists a Lipschitz family of real analytic symplectic coordinate transformations $\Phi_{v}:D_{v+1}\times\mathcal{O}_{v}^{*}\rightarrow D_{v}$
 and a closed subset
\begin{equation*}
 \mathcal{O}_{v+1}^{*}=\{(\omega^{b_{v+1}},\omega_{b_{v+1}}^{'}):\omega^{b_{v+1}}\in\mathcal{O}_{*}^{b_{v+1}}\}
\end{equation*}
of $\mathcal{O}_{v}^{*}$, where
 $$\mathcal{O}_{*}^{b_{v+1}}=\mathcal{O}^{v+1}\backslash\big(\mathop{\bigcup}\limits_{k,l}\mathcal{R}_{kl}^{v+1}\big),$$
 $$\mathcal{O}^{v+1}=\{\omega^{b_{v+1}}:(\omega^{b_{v+1}},\omega_{b_{v+1}}^{'})\in\mathcal{O}_{v}^{*}\}\subset[0,1]^{b_{v+1}},$$
$$\mathcal{R}^{v+1}_{kl}=\left\{\omega^{b_{v+1}}\in\mathcal{O}^{v+1}:|\langle k,\omega^{b_{v+1}}\rangle+\langle l,\Omega_{v+1}\rangle|< \frac{\alpha_{v+1}\langle l\rangle_{2}}{\big(1+(v+1)^2\big)(|k|+1)^{2b_{v+1}+2}}\right\},~~(k,l)\in\mathcal{Z}^{b_{v+1}},$$
such that for $H_{v+1}=H_{v}\circ\Phi_{v}=N_{v+1}+P_{v+1}$
the same assumptions are satisfied with $ v+1 $ in place of $v$.\\

\textbf{Remark 3.2} The assumption (\textbf{H2}) and Lemma 2.2 imply that $\tilde{P}^{b_{n}}$ in (\ref{Pn}) satisfy $~~~~$    $|X_{\tilde{P}^{b_{n}}}|_{r_{0},D(s_{0},r_{0})}^{\lambda_{0}}\leq C~~ (n\geq0)$. Here and later, the letter $C$ denotes suitable (possibly different) constants which are independent of iteration steps.

\subsection{Solving homological equations}

The coordinate transformation $\Phi_{v}$ is obtained as the time-1-map $X_{\mathcal{F}_{v}}^{t}|_{t=1}$ of the Hamiltonian vector field $X_{\mathcal{F}_{v}}$, where $\mathcal{F}_{v}$ has a similar expression of $R^{b_{v}}$,
\begin{equation*}
\begin{aligned}\label{}
&\mathcal{F}_{v}(\theta^{b_{v}},z,\bar{z},\omega^{b_{v}})
=\varepsilon_{v}F_{v}\\
=&\varepsilon_{v}F^{00b_{v}}+\varepsilon_{v}\langle F^{10b_{v}},z\rangle+\varepsilon_{v}\langle F^{01b_{v}},\bar{z}\rangle+\varepsilon_{v}\langle F^{20b_{v}}z,z\rangle+\varepsilon_{v}\langle F^{11b_{v}}z,\bar{z}\rangle+\varepsilon_{v}\langle F^{02b_{v}}\bar{z},\bar{z}\rangle\\
=&\varepsilon_{v}\sum_{0\neq k\in\mathbb{Z}^{b_{v}}}F_{k}^{00b_{v}}e^{\sqrt{-1}\langle k,\theta^{b_{v}}\rangle}+\varepsilon_{v}\sum_{j\geq0,k\in\mathbb{Z}^{b_{v}}}(F_{kj}^{10b_{v}}z_{j}+F_{kj}^{01b_{v}}\bar{z}_{j})e^{\sqrt{-1}\langle k,\theta^{b_{v}}\rangle}
\\
&+\varepsilon_{v}\sum_{i,j\geq0,k\in\mathbb{Z}^{b_{v}},|k|+|i-j|\neq0}(F_{kij}^{20b_{v}}z_{i}z_{j}+F_{kij}^{11b_{v}}z_{i}\bar{z}_{j}+F_{kij}^{02b_{v}}\bar{z}_{i}\bar{z}_{j})e^{\sqrt{-1}\langle k,\theta^{b_{v}}\rangle}.
\end{aligned}
\end{equation*}

\par By Taylor's formula, we have
\begin{equation}
\begin{aligned}\label{eqv}
H_{v+1}:=&H_{v}\circ \Phi_{v}\\
=&N_{v}+\varepsilon_{v}\{N_{v},F_{v}\}+\varepsilon_{v}^2\int_{0}^{1}(1-t)\{\{N_{v},F_{v}\},F_{v}\}\circ X_{\mathcal{F}_{v}}^{t}dt\\
&+\varepsilon_{v}R^{b_{v}}+\varepsilon_{v}^2\int_{0}^{1}\{R^{b_{v}},F_{v}\}\circ X_{\mathcal{F}_{v}}^{t}dt\\
&+\hat{P}_{v}^{high}+\varepsilon_{v}\{\hat{P}_{v}^{high},F_{v}\}+\varepsilon_{v}^2\int_{0}^{1}(1-t)\{\{\hat{P}_{v}^{high},F_{v}\},F_{v}\}\circ X_{\mathcal{F}_{v}}^{t}dt\\
&+(P_{v}-\hat{P}_{v})\circ \Phi_{v}.\\
\end{aligned}
\end{equation}
Then we obtain the modified homological equation

\begin{equation}\label{eq3.5}
 \begin{aligned}
\varepsilon_{v}\{N_{v},F_{v}\}+\varepsilon_{v}R^{b_{v}}+\varepsilon_{v}\{\hat{P}_{v}^{high},F_{v}\}^{low}=N_{v+1}-N_{v}.
 \end{aligned}
 \end{equation}
If the homological equation is solved, then the new perturbation term $P_{v+1}$ can be written as
\begin{eqnarray}
  P_{v+1}&=&\hat{P}_{v}^{high}+\varepsilon_{v}\{\hat{P}_{v}^{high},F_{v}\}^{high}\label{aa}\\
  &&+\varepsilon_{v}^2\int_{0}^{1}(1-t)\{\{N_{v}+\hat{P}_{v}^{high},F_{v}\},F_{v}\}\circ X_{\mathcal{F}_{v}}^{t}dt\label{bb}\\
  &&+\varepsilon_{v}^2\int_{0}^{1}\{R^{b_{v}},F_{v}\}\circ X_{\mathcal{F}_{v}}^{t}dt+(P_{v}-\hat{P}_{v})\circ \Phi_{v}.\label{cc}
\end{eqnarray}
Note that the terms in (\ref{aa}) have at least three normal variables. The terms in (\ref{aa}) will be left since they have no effect on the tori. To make the terms in (\ref{bb}) and (\ref{cc}) smaller with KAM iteration, different from the standard strategy in
J. P\"{o}schel \cite{Poschel1996}, we shrink the analytic radius of $z$ more slowly than \cite{Poschel1996} such that the final analytic radius of $z$ will be $r_{0}/2$ instead of $0$ (see the expression of the iterative constant $r_{v}$). Thus we can obtain a non-degenerate normal form in the end.

\par To solve the homological equation (\ref{eq3.5}), we should know the term $\{\hat{P}_{v}^{high},F_{v}\}^{low}$ exactly.\\
Let $\hat{P}_{v}^{high}=\hat{P}_{v0}^{high}+\hat{P}_{v1}^{high}$, where
\begin{equation*}\label{}
  \hat{P}_{v0}^{high}=\sum_{|\gamma|_{1}+|\kappa|_{1}=3}\hat{P}_{v}^{\gamma\kappa}(\theta^{b_{v}},\omega^{b_{v}})z^{\gamma}\bar{z}^{\kappa},~~\hat{P}_{v1}^{high}=\sum_{|\gamma|_{1}+|\kappa|_{1}\geq4}\hat{P}_{v}^{\gamma\kappa}(\theta^{b_{v}},\omega^{b_{v}})z^{\gamma}\bar{z}^{\kappa}.
\end{equation*}

Set $F_{v}=F_{v}^{0}+F_{v}^{1}+F_{v}^{2}$, where
\begin{eqnarray*}
  &&F_{v}^{0}=\sum_{0\neq k\in\mathbb{Z}^{b_{v}}}F_{k}^{00b_{v}}e^{\sqrt{-1}\langle k,\theta^{b_{v}}\rangle}, \\
  &&F_{v}^{1}=\sum_{j\geq0,k\in\mathbb{Z}^{b_{v}}}(F_{kj}^{10b_{v}}z_{j}+F_{kj}^{01b_{v}}\bar{z}_{j})e^{\sqrt{-1}\langle k,\theta^{b_{v}}\rangle}, \\
  &&F_{v}^{2}=\sum_{i,j\geq0,k\in\mathbb{Z}^{b_{v}},|k|+|i-j|\neq0}(F_{kij}^{20b_{v}}z_{i}z_{j}+F_{kij}^{11b_{v}}z_{i}\bar{z}_{j}+F_{kij}^{02b_{v}}\bar{z}_{i}\bar{z}_{j})e^{\sqrt{-1}\langle k,\theta^{b_{v}}\rangle}.
\end{eqnarray*}
 Denote $W^{b_{v}}=\{\hat{P}_{v}^{high},F_{v}\}^{low}$, then by a direct calculation, we obtain
\begin{equation*}
 W^{b_{v}}=\sqrt{-1}\sum_{j\geq0}\big(\partial_{z_{j}}\hat{P}_{v0}^{high}\partial_{\bar{z}_{j}}F_{v}^1-\partial_{\bar{z}_{j}}\hat{P}_{v0}^{high}\partial_{z_{j}}F_{v}^1\big)=\{\hat{P}_{v0}^{high},F_{v}^1\},
\end{equation*}
which is degree two in variables $(z,\bar{z}).$ Write $W^{b_{v}}=(W^{b_{v}})^0+(W^{b_{v}})^1+(W^{b_{v}})^2.$ Then we can easily get
\begin{equation}\label{}
\begin{aligned}
  W^{b_{v}}=(W^{b_{v}})^2&=\langle W^{20b_{v}}z,z\rangle+\langle W^{11b_{v}}z,\bar{z}\rangle+\langle W^{02b_{v}}\bar{z},\bar{z}\rangle.\\
\end{aligned}
\end{equation}
Let $(R^{b_{v}})^2$ be of the same form as $(W^{b_{v}})^2$ and set $B^{b_{v}}=(R^{b_{v}})^2+(W^{b_{v}})^2.$ More precisely,
\begin{equation*}
  B^{b_{v}}=\langle B^{20b_{v}}z,z\rangle+\langle B^{11b_{v}}z,\bar{z}\rangle+\langle B^{02b_{v}}\bar{z},\bar{z}\rangle
\end{equation*}
with
\begin{equation*}
  B^{20b_{v}}=R^{20b_{v}}+W^{20b_{v}},~~B^{11b_{v}}=R^{11b_{v}}+W^{11b_{v}},~~B^{02b_{v}}=R^{02b_{v}}+W^{02b_{v}}.
\end{equation*}

By the definition of $F_{v}$ and $b_{v}<b_{v}+1,$ it implies that $\{\mathop{\sum}\limits_{j\geq b_{v}+1}\omega_{i_{j}}J_{i_{j}},F_{v}\}=0.$
Moreover, it is easy to see that $\mathcal{F}_{v}$ and $\tilde{P}^{b_{n}}~(n\geq0)$ are independent of $J$. Therefore, (\ref{eq3.5}) is equivalent to the following homological equations:
$$
\left\{
     \begin{array}{ll}
      \sqrt{-1}\langle k,\omega^{b_{v}}\rangle F^{00b_{v}}_{k}=R^{00b_{v}}_{k},~~k\neq0,\\
      \sqrt{-1}\big(\langle k,\omega^{b_{v}}\rangle +\Omega_{jv}\big)F^{10b_{v}}_{kj}=R^{10b_{v}}_{kj},\\
      \sqrt{-1}\big(\langle k,\omega^{b_{v}}\rangle-\Omega_{jv}\big)F^{01b_{v}}_{kj}=R^{01b_{v}}_{kj},\\
     \sqrt{-1}\big(\langle k,\omega^{b_{v}}\rangle+\Omega_{iv}+\Omega_{jv}\big)F^{20b_{v}}_{kij}=R^{20b_{v}}_{kij}+W^{20b_{v}}_{kij}= B^{20b_{v}}_{kij},\\
     \sqrt{-1}\big(\langle k,\omega^{b_{v}}\rangle+\Omega_{iv}-\Omega_{jv}\big)F^{11b_{v}}_{kij}=R^{11b_{v}}_{kij}+W^{11b_{v}}_{kij}= B^{11b_{v}}_{kij},~~|k|+|i-j|\neq0,\\
     \sqrt{-1}\big(\langle k,\omega^{b_{v}}\rangle-\Omega_{iv}-\Omega_{jv}\big)F^{02b_{v}}_{kij}=R^{02b_{v}}_{kij}+W^{02b_{v}}_{kij}= B^{02b_{v}}_{kij},
     \end{array}
   \right.$$
with $\varepsilon_{v}\hat{N}_{v}\triangleq N_{v+1}-N_{v}=\varepsilon_{v}[R^{00b_{v}}]+\varepsilon_{v}\sum_{j\geq0}[B^{11b_{v}}_{jj}]z_{j}\bar{z}_{j}$, where
\begin{equation*}
  [R^{00b_{v}}]=\frac{1}{(2\pi)^{b_{v}}}\int_{\mathbb{T}^{b_{v}}}R^{00b_{v}}(\theta^{b_{v}},\omega^{b_{v}})d\theta^{b_{v}}
\end{equation*}
 will be omitted from $H_{v+1}$ in the following since it dose not affect the dynamics of the Hamiltonian vector field $X_{H_{v+1}}$ and $[B^{11b_{v}}_{jj}]$ is defined analogously.

Concerning the estimate of $X_{F_{v}}$, we have the following lemma.\\

\textbf{Lemma 3.3}\label{lem3.1}
Suppose that uniformly on $\mathcal{O}_{v}^{*}$,
$$|\langle k,\omega^{b_{v}}\rangle+\langle l,\Omega_{v}\rangle|\geq\frac{\alpha_{v}\langle l\rangle_{2}}{(1+v^2)(|k|+1)^{2b_{v}+2}},~~~~  (k,l)\in\mathcal{Z}^{b_{v}}.$$
 Then the linearized equation (\ref{eq3.5}) has a solution $F_{v}$ satisfying
\begin{equation*}
\begin{aligned}
|X_{F_{v}}|_{r_{v},D(s_{v}-3\sigma_{v},r_{v}-\rho_{v})\times\mathcal{O}_{v}^{*}}^{\lambda_{v}}\leq C\alpha_{v}^{-2}\sigma_{v}^{-1}A_{v}^2|X_{R^{b_{v}}}|^{\lambda_{v}}_{r_{v},D(s_{v},r_{v})\times\mathcal{O}_{v}^{*}}
\end{aligned}
\end{equation*}
with $A_{v}=\big(\frac{16(2b_{v}+3)}{e}\big)^{4b_{v}+6}\sigma_{v}^{-(5b_{v}+6)}.$\\

\textbf{Proof}  Consider the term $F^{10b_{v}},$ we note that $R^{10b_{v}}$ is an analytic map in $\ell^{a,p+2}$ with a Fourier series expansion whose coefficients $R_{k}^{10b_{v}}=(R_{k0}^{10b_{v}},R_{k1}^{10b_{v}},\cdots)$ satisfy
 \begin{equation}\label{a}
 \begin{aligned}
  \sum_{k\in\mathbb{Z}^{b_{v}}}\|R_{k}^{10b_{v}}\|^2_{a,p+2}e^{2|k|s_{v}}\leq2^{b_{v}}\|R^{10b_{v}}\|_{a,p+2,D(s_{v})}^2,\\
 \end{aligned}
 \end{equation}
 where $D(s_{v})=\{|\text{Im}\theta^{b_{v}}|<s_{v}\}$ and $\|R^{10b_{v}}\|_{a,p+2,D(s_{v})}=\mathop{\text{sup}}\limits_{D(s_{v})}\|R^{10b_{v}}\|_{a,p+2}.$
 According to the small divisor assumptions, we easily get
 \begin{equation*}\label{10}
 \begin{aligned}
  |F^{10b_{v}}_{kj}|\leq\alpha_{v}^{-1}(1+v^2)(|k|+1)^{2b_{v}+2}|R^{10b_{v}}_{kj}|
 \end{aligned}
 \end{equation*}
 and
\begin{equation}\label{b}
 \begin{aligned}
  \|F_{k}^{10b_{v}}\|_{a,p+2}^2&=|F_{k0}^{10b_{v}}|^2+\sum_{j\geq1}|F_{kj}^{10b_{v}}|^2j^{2(p+2)}e^{2ja}\\
  &\leq[\alpha_{v}^{-1}(1+v^2)(|k|+1)^{2b_{v}+2}]^2\|R_{k}^{10b_{v}}\|_{a,p+2}^2\\
 \end{aligned}
 \end{equation}
uniformly on $\mathcal{O}_{v}^{*}$. From (\ref{a}), (\ref{b}) and Lemma 7.3 in Appendix, it follows that
\begin{equation}\label{eq3.41}
 \begin{aligned}
  \|F^{10b_{v}}\|_{a,p+2,D(s_v-\sigma_v)}&\leq\sum_{k\in\mathbb{Z}^{b_{v}}}\|F_{k}^{10b_{v}}\|_{a,p+2}e^{|k|(s_{v}-\sigma_{v})}\\
  &\leq\sqrt{\sum_{k\in\mathbb{Z}^{b_{v}}}\|R_{k}^{10b_{v}}\|^2_{a,p+2}e^{2|k|s_{v}}}\sqrt{\sum_{k\in\mathbb{Z}^{b_{v}}}[\alpha_{v}^{-1}(1+v^2)(|k|+1)^{2b_{v}+2}]^2e^{-2|k|\sigma_{v}}}\\
  &\leq \alpha_{v}^{-1}A_{v}\|R^{10b_{v}}\|_{a,p+2,D(s_{v})}
 \end{aligned}
 \end{equation}
and
\begin{equation*}
  \|\partial_{\theta}F^{10b_{v}}\|_{a,p+2,D(s_v-\sigma_v)}\leq\alpha_{v}^{-1}A_{v}\|R^{10b_{v}}\|_{a,p+2,D(s_{v})}.
\end{equation*}
Since $R^{10b_{v}}=(R^{b_{v}})_{z}|_{z=\bar{z}=0}$, we can easily get $\|R^{10b_{v}}\|_{a,p+2,D(s_{v})}\leq r_{v}|X_{R^{b_{v}}}|_{r_{v},D(s_{v},r_{v})}$, and
 \begin{equation}\label{eql}
   |X_{\langle F^{10b_{v}},z\rangle}|_{r_{v},D(s_v-\sigma_v,r_{v})\times\mathcal{O}_{v}^{*}}\leq C\alpha_{v}^{-1}A_{v}|X_{R^{b_{v}}}|_{r_{v},D(s_{v},r_{v})\times\mathcal{O}_{v}^{*}}.
 \end{equation}

To estimate the Lipschitz semi-norm of $F^{10b_{v}},$ let $\delta_{kj}=\langle k,\omega^{b_{v}}\rangle+\Omega_{vj}$ and $\Delta=\Delta_{\omega^{b_{v}}\tilde{\omega}^{b_{v}}}$ for $\omega^{b_{v}},$ $\tilde{\omega}^{b_{v}}\in\mathcal{O}_{\ast}^{b_{v}}$, then we have

\begin{equation*}
\begin{aligned}
 \Delta{F}^{10b_{v}}_{kj}
 =-\frac{\sqrt{-1}\Delta{R}^{10b_{v}}_{kj}}{\delta_{kj}(\omega^{b_{v}})}+\frac{\sqrt{-1}{R}^{10b_{v}}_{kj}(\tilde{\omega}^{b_{v}})\Delta\delta_{kj}}{\delta_{kj}(\omega^{b_{v}})\delta_{kj}(\tilde{\omega}^{b_{v}})}.\label{3.11}
\end{aligned}
\end{equation*}
The small divisor assumptions imply that
\begin{equation*}
\begin{aligned}
|\Delta F^{10b_{v}}_{kj}|\leq&\big(\alpha_{v}^{-1}(1+v^2)(|k|+1)^{2b_{v}+2}\big)^2|{R}^{10b_{v}}_{kj}|(|k||\Delta\omega^{b_{v}}|+|\Delta\Omega_{vj}|j^{-2})\\
 &+\alpha_{v}^{-1}(1+v^2)(|k|+1)^{2b_{v}+2}|\Delta{R}^{10b_{v}}_{kj}|\label{3.11}
\end{aligned}
\end{equation*}
on $\mathcal{O}_{v}^{\ast}.$
Hence,
\begin{equation*}\label{9}
 \begin{aligned}
  \|\Delta F_{k}^{10b_{v}}\|_{a,p+2}\leq&\big(\alpha_{v}^{-1}(1+v^2)(|k|+1)^{2b_{v}+2}\big)^2\|{R}^{10b_{v}}_{k}\|_{a,p+2}(|k||\Delta\omega^{b_{v}}|+|\Delta\Omega_{v}|_{-\delta})\\
 &+\alpha_{v}^{-1}(1+v^2)(|k|+1)^{2b_{v}+2}\|\Delta{R}^{10b_{v}}_{k}\|_{a,p+2}.
 \end{aligned}
 \end{equation*}
Summing up the Fourier series as (\ref{eq3.41}), we have
\begin{equation*}\label{9}
 \begin{aligned}
  \|\Delta F^{10b_{v}}\|_{a,p+2,D(s_v-\sigma_v)}
  \leq& \alpha_{v}^{-2}A_{v}\|R^{10b_{v}}\|_{a,p+2,D(s_{v})}(|\Delta\omega^{b_{v}}|+|\Delta\Omega_{v}|_{-\delta})\\
  &+\alpha_{v}^{-1}A_{v}\|\Delta R^{10b_{v}}\|_{a,p+2,D(s_{v})}.\\
 \end{aligned}
 \end{equation*}
Dividing by $|\omega^{b_{v}}-\tilde{\omega}^{b_{v}}|$ and taking the supremum over $\omega^{b_{v}}\neq\tilde{\omega}^{b_{v}}$ in $\mathcal{O}_{*}^{b_{v}}$ , we obtain
\begin{equation*}\label{eq3.42}
\begin{aligned}
\|F^{10b_{v}}\|^{lip}_{a,p+2,D(s_{v}-\sigma_{v})}&\leq \alpha_{v}^{-1}A_{v}\Big(\frac{M_{v}+1}{\alpha_{v}}\|R^{10b_{v}}\|_{a,p+2,D(s_{v})}
+\|{R}^{10b_{v}}\|^{lip}_{a,p+2,D(s_{v})}\Big)
\end{aligned}
\end{equation*}
and
\begin{equation*}\label{eq3.42}
\begin{aligned}
\|\partial_{\theta}F^{10b_{v}}\|^{lip}_{a,p+2,D(s_{v}-\sigma_{v})}&\leq \alpha_{v}^{-1}A_{v}\Big(\frac{M_{v}+1}{\alpha_{v}}\|R^{10b_{v}}\|_{a,p+2,D(s_{v})}
+\|{R}^{10b_{v}}\|^{lip}_{a,p+2,D(s_{v})}\Big)
\end{aligned}
\end{equation*}
in view of (\ref{eq3.1}). Then, we have
\begin{equation*}
\begin{aligned}\label{eqd}
|X_{\langle F^{10b_{v}},z\rangle}|_{r_{v},D(s_{v}-\sigma_{v},r_{v})\times\mathcal{O}_{v}^{*}}^{lip}
\leq C\alpha_{v}^{-1}A_{v}\Big(\frac{M_{v}+1}{\alpha_{v}}|X_{R^{b_{v}}}|_{r_{v},D(s_{v},r_{v})\times\mathcal{O}_{v}^{*}}+|X_{R^{b_{v}}}|_{r_{v},D(s_{v},r_{v})\times\mathcal{O}_{v}^{*}}^{lip}\Big).
\end{aligned}
\end{equation*}
Thus, together with (\ref{eql}), we arrive at
\begin{equation}
\begin{aligned}\label{F1}
|X_{\langle F^{10b_{v}},z\rangle}|_{r_{v},D(s_{v}-\sigma_{v},r_{v})\times\mathcal{O}_{v}^{*}}^{\lambda_{v}}
\leq C\alpha_{v}^{-1}A_{v}|X_{R^{b_{v}}}|^{\lambda_{v}}_{r_{v},D(s_{v},r_{v})\times\mathcal{O}_{v}^{*}}.
\end{aligned}
\end{equation}
For the term $F^{01b_{v}},$ the same estimate as (\ref{F1}) can be obtained, thus, we get
\begin{equation}\label{eqa}
   |X_{F_{v}^{1}}|^{\lambda_{v}}_{r_{v},D(s_v-\sigma_v,r_{v})\times\mathcal{O}_{v}^{*}}\leq C\alpha_{v}^{-1}A_{v}|X_{R^{b_{v}}}|^{\lambda_{v}}_{r_{v},D(s_{v},r_{v})\times\mathcal{O}_{v}^{*}}.
 \end{equation}

Before considering the term $F^{11b_{v}}$, we should get the estimate of the term $B^{b_{v}}$.
Recall that $$B^{b_{v}}=(R^{b_{v}})^{2}+(W^{b_{v}})^{2},$$ and $$(W^{b_{v}})^{2}=\{\hat{P}_{v}^{high},F_{v}\}^{low}=\{\hat{P}_{v0}^{high},F_{v}^{1}\}.$$
Then by (\ref{Ph}), (\ref{eqi}), (\ref{eqa}) and the generalized Cauchy inequality, we obtain

\begin{equation}
\begin{aligned}\label{Bv}
|X_{B^{b_{v}}}|_{r_{v},D(s_{v}-2\sigma_{v},r_{v}-d_{v})\times\mathcal{O}_{v}^{*}}^{\lambda_{v}}&\leq |X_{(R^{b_{v}})^2}|_{r_{v},D(s_{v},r_{v})\times\mathcal{O}_{v}^{*}}^{\lambda_{v}}+|X_{(W^{b_{v}})^2}|_{r_{v},D(s_{v}-2\sigma_{v},r_{v}-d_{v})\times\mathcal{O}_{v}^{*}}^{\lambda_{v}}\\
&\leq C|X_{R^{b_{v}}}|_{r_{v},D(s_{v},r_{v})\times\mathcal{O}_{v}^{*}}^{\lambda_{v}}+C\sigma_{v}^{-1}|X_{\hat{P}_{v}^{high}}|_{r_{v},D(s_{v},r_{v})\times\mathcal{O}_{v}^{*}}^{\lambda_{v}}|X_{F_{v}^{1}}|_{r_{v},D(s_{v}-\sigma_{v},r_{v})\times\mathcal{O}_{v}^{*}}^{\lambda_{v}}\\
&\leq C|X_{R^{b_{v}}}|_{r_{v},D(s_{v},r_{v})\times\mathcal{O}_{v}^{*}}^{\lambda_{v}}\Big(1+\alpha_{v}^{-1}\sigma_{v}^{-1}A_{v}|X_{\hat{P}_{v}^{high}}|_{r_{v},D(s_{v},r_{v})\times\mathcal{O}_{v}^{*}}^{\lambda_{v}}\Big)\\
&\leq C\alpha_{v}^{-1}\sigma_{v}^{-1}A_{v}|X_{R^{b{v}}}|_{r_{v},D(s_{v},r_{v})\times\mathcal{O}_{v}^{*}}^{\lambda_{v}}.
\end{aligned}
\end{equation}
By the generalized Cauchy inequality, we have
$$\|B^{11b_{v}}\|_{a,p+2,p,D(s_{v}-2\sigma_{v})}\leq\frac{1}{r_{v}}\|(B^{b_{v}})_{z}\|_{a,p+2,D(s_v-2\sigma_{v},r_{v}-d_{v})}\leq|X_{B^{b_{v}}}|_{r_{v},D(s_{v}-2\sigma_{v},r_{v}-d_{v})\times\mathcal{O}_{v}^{*}},$$ where
$\|\cdot\|_{a,p+2,p,D(s_{v}-2\sigma_{v})}$ is the operator norm of bounded linear operators from $\ell^{a,p}$ to $\ell^{a,p+2}$.
 This is equivalent to that $\tilde{B}^{b_{v}}=(v_{i}B_{ij}^{11b_{v}}w_{j})$ is a bounded linear operator of $\ell^2$ into itself with the operator norm $\||\tilde{B}^{b_{v}}\||_{D(s_{v}-2\sigma_{v})}=\|B^{11b_{v}}\|_{a,p+2,p,D(s_{v}-2\sigma_{v})},$ where $v_{i}, w_{j}$ are certain weights (see \cite{Poschel1996}).
Expanding $\tilde{B}^{b_{v}}$ into its Fourier series and as before, we know that $\mathop{\sum}\limits_{k\in\mathbb{Z}^{b_{v}}}\|| \tilde{B}_{k}^{b_{v}}\||^2e^{2|k|(s_{v}-2\sigma_{v})}\leq2^{b_{v}}\||\tilde{B}^{b_{v}}\||_{D(s_{v}-2\sigma_{v})}$.
By the small divisor assumptions and $|i^{2}-j^2|=|i-j|(i+j)$, we find that the corresponding coefficient $\tilde{F}_{k}^{b_{v}}=(\tilde{F}_{kij}^{b_{v}})$ satisfies the following estimate
 \begin{equation*}\label{9}
 \begin{aligned}
  |\tilde{F}^{b_{v}}_{kij}|\leq\frac{\alpha_{v}^{-1}(1+v^2)(|k|+1)^{2b_{v}+2}}{|i-j|}|\tilde{B}^{b_{v}}_{kij}|,~~~~|k|+|i-j|\neq0,\\
 \end{aligned}
 \end{equation*}
 while $\tilde{B}^{b_{v}}_{0jj}=0.$ Lemma 7.2 in Appendix implies that
 $$|||\tilde{F}_{k}^{b_{v}}|||\leq3\alpha_{v}^{-1}(1+v^2)(|k|+1)^{2b_{v}+2}|||\tilde{B}_{k}^{b_{v}}|||$$
 uniformly on $\mathcal{O}_{v}^{\ast}.$
 Summing up the Fourier series as before,
 \begin{equation*}\label{10}
\begin{aligned}
 \||\tilde{F}^{b_{v}}\||_{D(s_v-3\sigma_v)}\leq\sum_{k\in\mathbb{Z}^{b_{v}}}\||\tilde{F}_{k}^{b_{v}}\||e^{|k|(s_v-3\sigma_v)}\leq 3\alpha_{v}^{-1}A_{v}\||\tilde{B}^{b_{v}}\||_{D(s_v-2\sigma_{v})}\\
\end{aligned}
\end{equation*}
and
 \begin{equation*}\label{}
\begin{aligned}
 \||\partial_{\theta}\tilde{F}^{b_{v}}\||_{D(s_v-3\sigma_v)}\leq\sum_{k\in\mathbb{Z}^{b_{v}}}|k|\||\tilde{F}_{k}^{b_{v}}\||e^{|k|(s_v-3\sigma_v)}\leq 3\alpha_{v}^{-1}A_{v}\||\tilde{B}^{b_{v}}\||_{D(s_v-2\sigma_{v})}.\\
\end{aligned}
\end{equation*}
Thus,
 \begin{equation*}\label{10}
\begin{aligned}
 \|F^{11b_{v}}\|_{a,p+2,p,D(s_{v}-3\sigma_{v})}
 &\leq 3\alpha_{v}^{-1}A_{v}|X_{B^{b_{v}}}|_{r_{v},D(s_{v}-2\sigma_{v},r_{v}-d_{v})}
\end{aligned}
\end{equation*}
and
\begin{equation*}\label{10}
\begin{aligned}
 \|\partial_{\theta}F^{11b_{v}}\|_{a,p+2,p,D(s_{v}-3\sigma_{v})}
 &\leq 3\alpha_{v}^{-1}A_{v}|X_{B^{b_{v}}}|_{r_{v},D(s_{v}-2\sigma_{v},r_{v}-d_{v})}.
\end{aligned}
\end{equation*}
Finally, we have
\begin{equation}\label{eqb}
|X_{\langle F^{11b_{v}}z,\bar{z}\rangle}|_{r_{v},D(s_{v}-3\sigma_{v},r_{v}-d_{v})\times\mathcal{O}_{v}^{*}}\leq C\alpha_{v}^{-1}A_{v}|X_{B^{b_{v}}}|_{r_{v},D(s_{v}-2\sigma_{v},r_{v}-d_{v})\times\mathcal{O}_{v}^{*}}.
\end{equation}

\par To obtain the estimate of Lipschitz semi-norm, let $\delta_{kij}=\langle k,\omega\rangle+\Omega_{vi}-\Omega_{vj},$ then we have
\begin{equation*}
\begin{aligned}
 \Delta{\tilde{F}}^{b_{v}}_{kij}
 =-\frac{\sqrt{-1}\Delta{\tilde{B}}^{b_{v}}_{kij}}{\delta_{kij}(\omega^{b_{v}})}+\frac{\sqrt{-1}{\tilde{B}}^{b_{v}}_{kij}(\tilde{\omega}^{b_{v}})\Delta\delta_{kij}}{\delta_{kij}(\omega^{b_{v}})\delta_{kij}(\tilde{\omega}^{b_{v}})},\label{3.11}
\end{aligned}
\end{equation*}
which implies that
\begin{equation*}
\begin{aligned}
 |\Delta{\tilde{F}}^{b_{v}}_{kij}|\leq&\big(\alpha_{v}^{-1}(1+v^2)(|k|+1)^{2b_{v}+2}\big)^2\frac{|{\tilde{B}}^{b_{v}}_{kij}|}{|i-j|}\big(|k||\Delta\omega|+2|\Delta\Omega_{v}|_{-\delta}\big)\\
 &+\alpha_{v}^{-1}(1+v^2)(|k|+1)^{2b_{v}+2}\frac{|\Delta{\tilde{B}}^{b_{v}}_{kij}|}{|i-j|}.\label{3.11}
\end{aligned}
\end{equation*}
Thus, we obtain
\begin{equation*}
\begin{aligned}
 \||\Delta{\tilde{F}}^{b_{v}}\||_{D(s_{v}-3\sigma_{v})}&\leq 3\alpha_{v}^{-2}A_{v}\||{\tilde{B}}^{b_{v}}\||_{D(s_{v}-2\sigma_{v})}\big(|\Delta\omega|+2|\Delta\Omega_{v}|_{-\delta}\big)+3\alpha_{v}^{-1}A_{v}\||\Delta{\tilde{B}}^{b_{v}}\||_{D(s_{v}-2\sigma_{v})}.\label{3.11}
\end{aligned}
\end{equation*}
Dividing by $|\omega^{b_{v}}-\tilde{\omega}^{b_{v}}|$ and taking the supremum over $\omega^{b_{v}}\neq\tilde{\omega}^{b_{v}}$ in $\mathcal{O}_{\ast}^{b_{v}}$, we arrive at
\begin{equation*}
\begin{aligned}\label{eq3.11}
||F^{11b_{v}}\|^{lip}_{a,p+2,p,D(s_{v}-3\sigma_{v})}\leq 6\alpha_{v}^{-1}A_{v}\Big(\frac{M_{v}+1}{\alpha_{v}}\|B^{11b_{v}}\|_{a,p+2,p,D(s_{v}-2\sigma_{v})}
+\|B^{11b_{v}}\|^{lip}_{a,p+2,p,D(s_{v}-2\sigma_{v})}\Big).
\end{aligned}
\end{equation*}
In the same way as (\ref{eqb}), we have
\begin{equation*}
\begin{aligned}\label{eqc}
|X_{\langle F^{11b_{v}}z,\bar{z}\rangle}|_{r_{v},D(s_{v}-3\sigma_{v},r_{v}-d_{v})\times\mathcal{O}_{v}^{*}}^{lip}\leq& C\alpha_{v}^{-1}A_{v}\Big(\frac{M_{v}+1}{\alpha_{v}}|X_{B^{b_{v}}}|_{r_{v},D(s_{v}-2\sigma_{v},r_{v}-d_{v})\times\mathcal{O}_{v}^{*}}\\
&+|X_{B^{b_{v}}}|^{lip}_{r_{v},D(s_{v}-2\sigma_{v},r_{v}-d_{v})\times\mathcal{O}_{v}^{*}}\Big),
\end{aligned}
\end{equation*}
which, together with (\ref{Bv}) and (\ref{eqb}), leads to

\begin{equation}
\begin{aligned}\label{eq4.45}
|X_{\langle F^{11b_{v}}z,\bar{z}\rangle}|_{r_{v},D(s_{v}-3\sigma_{v},r_{v}-d_{v})\times\mathcal{O}_{v}^{*}}^{\lambda_{v}}&\leq C\alpha_{v}^{-1}A_{v}|X_{B^{b_{v}}}|^{\lambda_{v}}_{r_{v},D(s_{v}-2\sigma_{v},r_{v}-d_{v})\times\mathcal{O}_{v}^{*}}\\
&\leq C\alpha_{v}^{-2}\sigma_{v}^{-1}A_{v}^2|X_{R^{b{v}}}|_{r_{v},D(s_{v},r_{v})\times\mathcal{O}_{v}^{*}}^{\lambda_{v}}.
\end{aligned}
\end{equation}

For the other terms of $F_{v}$, the same estimates or even better ones than (\ref{F1}) and (\ref{eq4.45}) can be obtained. Thus, we finally get the estimate of the Hamiltonian vector field $X_{F_{v}}$
\begin{equation}
\begin{aligned}\label{eqe}
|X_{F_{v}}|_{r_{v},D(s_{v}-3\sigma_{v},r_{v}-d_{v})\times\mathcal{O}_{v}^{*}}^{\lambda_{v}}\leq C\alpha_{v}^{-2}\sigma_{v}^{-1}A_{v}^2|X_{R^{b{v}}}|_{r_{v},D(s_{v},r_{v})\times\mathcal{O}_{v}^{*}}^{\lambda_{v}}.
\end{aligned}
\end{equation}

The proof of Lemma 3.3 is completed.~~~~~$\Box$

Note that $|X_{R^{b_{v}}}|_{r_{v},D(s_{v},r_{v})\times\mathcal{O}_{v}^{*}}^{\lambda_{v}}\leq \frac{1}{2}.$
Then using the generalized Cauchy inequality, we get
\begin{equation}
\begin{aligned}\label{eq3.37}
\frac{1}{\sigma_{v}}|X_{\mathcal{F}_{v}}|^{\lambda_{v}}_{r_{v},D(s_{v}-3\sigma_{v},r_{v}-d_{v})\times\mathcal{O}_{v}^{*}},|DX_{\mathcal{F}_{v}}|^{\lambda_{v}}_{r_{v},r_{v},D(s_{v}-4\sigma_{v},r_{v}-2d_{v})\times\mathcal{O}_{v}^{*}}
\leq C \varepsilon_{v}^{1-\frac{1}{4}\rho}
\end{aligned}
\end{equation}
by $\varepsilon_{v}^{\frac{1}{8}\rho}\big(\sigma_{v}^{-1}A_{v}\big)^2|X_{R^{b_{v}}}|^{\lambda_{v}}_{r_{v},D(s_{v},r_{v})\times\mathcal{O}_{v}^{*}}\leq C$ as $\varepsilon\ll1$, where we require $d_{v}/r_{v}\geq\sigma_{v}$ which is fulfilled by setting $\sigma_{0}\leq1/20$, moreover $C$ is an absolute constant independent of $v$ and $\varepsilon$. Here we use the operator norm
\begin{equation}\label{l}
|L|_{r,s}=\mathop{\text{sup}}\limits_{W\neq0}\frac{|LW|_{r,a,p+2}}{|W|_{s,a,p}}
\end{equation}
with $|\cdot|_{r,a,p+2}$ defined in (\ref{eq3.35}), and $|\cdot|_{s,a,p}$ defined analogously.
Then the flow $X_{\mathcal{F}_{v}}^{t}$ of the vector field $X_{\mathcal{F}_{v}}$ exists on $D(s_{v}-4\sigma_{v},r_{v}-2d_{v})$ for $-1\leq t\leq1$ and takes this domain into $D(s_{v}-3\sigma_{v},r_{v}-d_{v})$. By Lemma A.4 in \cite{Poschel1996} and (\ref{eq3.37}), we obtain
\begin{equation}\label{eq3.33}
|X_{\mathcal{F}_{v}}^t-id|^{\lambda_{v}}_{r_{v},D(s_{v}-4\sigma_{v},r_{v}-2d_{v})\times\mathcal{O}_{v}^{*}}\leq C|X_{\mathcal{F}_{v}}|^{\lambda_{v}}_{r_{v},D(s_{v}-3\sigma_{v},r_{v}-d_{v})\times\mathcal{O}_{v}^{*}}\leq C\varepsilon_{v}^{1-\frac{1}{4}\rho}
\end{equation}
for $-1\leq t\leq 1$.
Similarly, the flow takes $D(s_{v}-5\sigma_{v},r_{v}-3d_{v})$ into $D(s_{v}-4\sigma_{v},r_{v}-2d_{v})$
and by generalized Cauchy inequality, we also have
\begin{equation}\label{eq3.34}
|DX_{\mathcal{F}_{v}}^{t}-Id|^{\lambda_{v}}_{r_{v},r_{v},D(s_{v}-5\sigma_{v},r_{v}-3d_{v})\times\mathcal{O}_{v}^{*}}<C\frac{1}{\sigma_{v}}|X_{\mathcal{F}_{v}}|^{\lambda_{v}}_{r_{v},D(s_{v}-3\sigma_{v},r_{v}-d_{v})\times\mathcal{O}_{v}^{*}}
\leq C\varepsilon_{v}^{1-\frac{1}{4}\rho}
\end{equation}
for $-1\leq t\leq1.$
\subsection{The new Hamiltonian}

From (\ref{eqv}) and (\ref{eq3.5}) we get the new Hamiltonian $H_{v}\circ\Phi_{v}=N_{v+1}+P_{v+1}$
with
\begin{equation}\label{eqj}
  N_{v+1}=N_{v}+\varepsilon_{v}\sum_{j\geq0}[B_{jj}^{11b_{v}}]z_{j}\bar{z}_{j}:=N_{v}+\varepsilon_{v}\langle\hat{\Omega}_{v},z\bar{z}\rangle
\end{equation}
 and
\begin{equation*}
\begin{aligned}\label{v}
P_{v+1}=&\varepsilon_{v}^{2}\int_{0}^{1}(1-t)\{\{N_{v},F_{v}\},F_{v}\}\circ X_{\mathcal{F}_{v}}^{t}dt\\
&+\varepsilon_{v}^{2}\int_{0}^{1}\{R^{b_{v}},F_{v}\}\circ X_{\mathcal{F}_{v}}^{t}dt+\varepsilon_{v}^2\int_{0}^{1}(1-t)\{\{\hat{P}_{v}^{high},F_{v}\},F_{v}\}\circ X_{\mathcal{F}_{v}}^{t}dt+(P_{v}-\hat{P}_{v})\circ \Phi_{v}\\
&+\hat{P}_{v}^{high}+\varepsilon_{v}\{\hat{P}_{v}^{high},F_{v}\}^{high}.\\
\end{aligned}
\end{equation*}
Denote $$\Omega_{(v+1)j}(\omega^{b_{v}})=\Omega_{vj}(\omega^{b_{v-1}})+\varepsilon_{v}[B_{jj}^{11b_{v}}](\omega^{b_{v}}).$$
Thus, we obtain the new normal form $$N_{v+1}=\langle\omega^{b_{v+1}},J^{b_{v+1}}\rangle+\sum_{n\geq v+2}\langle\omega_{1}^{b_{n}},J_{1}^{b_{n}}\rangle+\langle\Omega_{v+1},z\bar{z}\rangle.$$

Let

\begin{equation}\label{eqw}
  P_{v+1}=\hat{P}_{v+1}+\sum_{n\geq v+2}\varepsilon^{(1+\rho)^{n}}\tilde{P}^{b_{n}}\circ \Phi_{0} \circ \Phi_{1} \circ \cdots \circ \Phi_{v}
\end{equation}
with

\begin{equation}\label{eqx}
\hat{P}_{v+1}=\varepsilon_{v+1}R^{b_{v+1}}+\hat{P}_{v+1}^{high},
\end{equation}
\begin{equation}
\begin{aligned}\label{eqD}
\hat{P}_{v+1}^{high}=&\hat{P}_{v}^{high}+\varepsilon_{v}\{\hat{P}_{v}^{high},F_{v}\}^{high}
+\varepsilon_{v}^2\Big(\int_{0}^{1}(1-t)\{\{\hat{P}_{v}^{high},F_{v}\},F_{v}\}\circ X_{\mathcal{F}_{v}}^{t}dt\Big)^{high}\\
&+\varepsilon^{(1+\rho)^{v+1}}\Big(\tilde{P}^{b_{v+1}}\circ \Phi_{0} \circ \Phi_{1} \circ \cdots \circ \Phi_{v}\Big)^{high},
\end{aligned}
\end{equation}

and

\begin{equation*}\label{}
\begin{aligned}
  \varepsilon_{v+1}R^{b_{v+1}}
  =&\varepsilon_{v}^{2}\int_{0}^{1}(1-t)\{\{N_{v},F_{v}\},F_{v}\}\circ X_{\mathcal{F}_{v}}^{t}dt+\varepsilon_{v}^{2}\int_{0}^{1}\{R^{b_{v}},F_{v}\}\circ X_{\mathcal{F}_{v}}^{t}dt\\
  &+\varepsilon_{v}^{2}\Big(\int_{0}^{1}(1-t)\{\{\hat{P}_{v}^{high},F_{v}\},F_{v}\}\circ X_{\mathcal{F}_{v}}^{t}dt\Big)^{low}\\
  &+\varepsilon^{(1+\rho)^{v+1}}\Big(\tilde{P}^{b_{v+1}}\circ \Phi_{0} \circ \Phi_{1} \circ \cdots \circ \Phi_{v}\Big)^{low}.
\end{aligned}
\end{equation*}
Since $\hat{P}_{v}$ and $F_{v}$ depend on $z,$ $\bar{z}$, $\theta^{b_{v}}$ and $\omega^{b_{v}}$. Moreover, $\tilde{P}^{b_{v+1}}$ depends on $z$, $\bar{z}$, $\theta_{1}^{b_{v+1}}$ and $\omega_{1}^{b_{v+1}}$. It is easy to check that $\hat{P}_{v+1}$ depends on $z$, $\bar{z}$, $\theta^{b_{v+1}}$ and $\omega^{b_{v+1}}$. Thus, in order to remove the non-normalized terms in $(\hat{P}_{v+1})^{low}$ in the next step, we only need to treat finite frequencies.

\subsection{Estimate of the new norm form and new perturbation}
The aim of this section is to estimate the new normal form $N_{v+1}$ and the new perturbation $P_{v+1}$ in (\ref{eqw}).
Now we consider $R^{b_{v+1}}$ and we will prove that (\ref{eqi}) is fulfilled with $v+1$ in place of $v$. Concerning $R^{b_{v+1}},$ we have
\begin{equation}\label{eq3.13}
\begin{aligned}
X_{R^{b_{v+1}}}
=&\varepsilon_{v}^{1-\frac{1}{2}\rho}\int_{0}^{1}(1-t)(X_{\mathcal{F}_{v}}^t)^{*}[X_{\{N_{v},F_{v}\}},X_{F_{v}}]dt
+\varepsilon_{v}^{1-\frac{1}{2}\rho}\int_{0}^{1}(X_{\mathcal{F}_{v}}^t)^{*}[X_{R^{b_{v}}},X_{F_{v}}]dt\\
&+\varepsilon_{v}^{1-\frac{1}{2}\rho}X_{\big(\int_{0}^{1}(1-t)\{\{\hat{P}_{v}^{high},F_{v}\},F_{v}\}\circ X_{\mathcal{F}_{v}}^{t}dt\big)^{low}}
+\varepsilon_{v+1}^{-1}\varepsilon^{(1+\rho)^{v+1}}X_{\big(\tilde{P}^{b_{v+1}}\circ \Phi_{0} \circ \Phi_{1} \circ \cdots \circ \Phi_{v}\big)^{low}}.
\end{aligned}
\end{equation}

In the following, we will estimate every part of $X_{R^{b_{v+1}}}$, and we wish to show that $R^{b_{v+1}}$ satisfies (\ref{eqi}) with $v+1$ in place of $v$.
By the generalized Cauchy estimate, we obtain
\begin{equation}\label{eqy}
  |X_{\{N_{v},F_{v}\}}|^{\lambda_{v}}_{r_{v},D(s_{v}-4\sigma_{v},r_{v}-2d_{v})\times\mathcal{O}_{v}^{*}}\leq C\sigma_{v}^{-1}|X_{F_{v}}|^{\lambda_{v}}_{r_{v},D(s_{v}-3\sigma_{v},r_{v}-d_{v})\times\mathcal{O}_{v}^{*}}
\end{equation}

Following the same lines as (12) in \cite{Poschel1996}, we obtain that for any vector field $Y$,
\begin{equation}\label{eq3.17}
 |(X_{\mathcal{F}_{v}}^t)^*Y|^{\lambda_{v}}_{r_{v},D(s_{v}-6\sigma_{v},r_{v}-4d_{v})\times\mathcal{O}_{v}^{*}}<2|Y|^{\lambda_{v}}_{r_{v},D(s_{v}-5\sigma_{v},r_{v}-3d_{v})\times\mathcal{O}_{v}^{*}}~~~~(0\leq t\leq 1).
\end{equation}

Therefore, by  the generalized Cauchy inequality and (\ref{eqi}), (\ref{eqy}), we get

\begin{equation}\label{eqz}
\begin{aligned}
  &|(X_{\mathcal{F}_{v}}^t)^{*}[X_{\{N_{v},F_{v}\}},X_{F_{v}}]|^{\lambda_{v}}_{r_{v},D(s_{v}-6\sigma_{v},r_{v}-4d_{v})\times\mathcal{O}_{v}^{*}}\\
  \leq& C\sigma_{v}^{-1}|X_{\{N_{v},F_{v}\}}|^{\lambda_{v}}_{r_{v},D(s_{v}-4\sigma_{v},r_{v}-2d_{v})\times\mathcal{O}_{v}^{*}}|X_{F_{v}}|^{\lambda_{v}}_{r_{v},D(s_{v}-3\sigma_{v},r_{v}-d_{v})\times\mathcal{O}_{v}^{*}}\\
  \leq& C\sigma_{v}^{-2}\Big(|X_{F_{v}}|^{\lambda_{v}}_{r_{v},D(s_{v}-3\sigma_{v},r_{v}-d_{v})\times\mathcal{O}_{v}^{*}}\Big)^2
\end{aligned}
\end{equation}
and
\begin{equation}\label{eq3.7}
\begin{aligned}
|(X_{\mathcal{F}_{v}}^t)^{*}[X_{R^{b_{v}}},X_{F_{v}}]|^{\lambda_{v}}_{r_{v},D(s_{v}-6\sigma_{v},r_{v}-4d_{v})\times\mathcal{O}_{v}^{*}}
&\leq2|[X_{R^{b_{v}}},X_{F_{v}}]|^{\lambda_{v}}_{r_{v},D(s_{v}-5\sigma_{v},r_{v}-3d_{v})\times\mathcal{O}_{v}^{*}}\\
&\leq C\sigma_{v}^{-1}|X_{F_{v}}|^{\lambda_{v}}_{r_{v},D(s_{v}-3\sigma_{v},r_{v}-d_{v})\times\mathcal{O}_{v}^{*}}.\\
\end{aligned}
\end{equation}
In the same way as (\ref{eqz}) (\ref{eq3.7}) and by (\ref{Ph}), we have
\begin{equation}\label{eqA}
\begin{aligned}
&|X_{\int_{0}^{1}(1-t)\{\{\hat{P}_{v}^{high},F_{v}\},F_{v}\}\circ X_{\mathcal{F}_{v}}^{t}dt}|^{\lambda_{v}}_{r_{v},D(s_{v}-6\sigma_{v},r_{v}-4d_{v})\times\mathcal{O}_{v}^{*}}\\
\leq&|(X_{\mathcal{F}_{v}}^t)^{*}[X_{\{\hat{P}_{v}^{high},F_{v}\}},X_{F_{v}}]|^{\lambda_{v}}_{r_{v},D(s_{v}-6\sigma_{v},r_{v}-4d_{v})\times\mathcal{O}_{v}^{*}}\\
\leq& C\sigma_{v}^{-1}|X_{\{\hat{P}_{v}^{high},F_{v}\}}|^{\lambda_{v}}_{r_{v},D(s_{v}-4\sigma_{v},r_{v}-2d_{v})\times\mathcal{O}_{v}^{*}}|X_{F_{v}}|^{\lambda_{v}}_{r_{v},D(s_{v}-3\sigma_{v},r_{v}-d_{v})\times\mathcal{O}_{v}^{*}}\\
\leq& C\sigma_{v}^{-2}\Big(|X_{F_{v}}|^{\lambda_{v}}_{r_{v},D(s_{v}-3\sigma_{v},r_{v}-d_{v})\times\mathcal{O}_{v}^{*}}\Big)^2.
\end{aligned}
\end{equation}

Then (\ref{eqA}) implies that
\begin{equation}\label{eqB}
\begin{aligned}
\Big|X_{\big(\int_{0}^{1}(1-t)\{\{\hat{P}_{v}^{high},F_{v}\},F_{v}\}\circ X_{\mathcal{F}_{v}}^{t}dt\big)^{low}}\Big|^{\lambda_{v}}_{r_{v},D(s_{v}-6\sigma_{v},r_{v}-4d_{v})\times\mathcal{O}_{v}^{*}}\leq C\sigma_{v}^{-2}\Big(|X_{F_{v}}|^{\lambda_{v}}_{r_{v},D(s_{v}-3\sigma_{v},r_{v}-d_{v})\times\mathcal{O}_{v}^{*}}\Big)^2.
\end{aligned}
\end{equation}

By repeatedly applying (\ref{eq3.17}) to $\tilde{P}^{b_{v+1}}\circ \Phi_{0} \circ \Phi_{1} \circ \cdots \circ \Phi_{v}$, we have

\begin{equation}\label{eqg}
\begin{aligned}
  |X_{\tilde{P}^{b_{v+1}}\circ \Phi_{0} \circ \Phi_{1} \circ \cdots \circ \Phi_{v}}|^{\lambda_{v}}_{r_{v},D(s_{v}-6\sigma_{v},r_{v}-4d_{v})\times\mathcal{O}_{v}^{*}}
   &\leq2^{v+3}C|X_{\tilde{P}^{b_{v+1}}}|^{\lambda_{0}}_{r_{0},D(s_{0},r_{0})\times\mathcal{O}_{v}^{*}}.
  \end{aligned}
\end{equation}
Thus (\ref{eqg}) and Remark 3.2  imply that

\begin{equation}\label{eqC}
\begin{aligned}
  \Big|X_{\big(\tilde{P}^{b_{v+1}}\circ \Phi_{0} \circ \Phi_{1} \circ \cdots \circ \Phi_{v}\big)^{low}}\Big|^{\lambda_{v}}_{r_{v},D(s_{v}-6\sigma_{v},r_{v}-4d_{v})\times\mathcal{O}_{v}^{*}}
   &\leq2^{v+3}C.
  \end{aligned}
\end{equation}

So altogether we obtain
\begin{equation}\label{eqh}
\begin{aligned}
|X_{R^{b_{v+1}}}|^{\lambda_{v}}_{r_{v},D(s_{v}-6\sigma_{v}, r_{v}-4d_{v})\times\mathcal{O}_{v}^{*}}
\leq& C\varepsilon_{v}^{1-\frac{1}{2}\rho}\sigma_{v}^{-2}\big(\alpha_{v}^{-2}\sigma_{v}^{-1}A_{v}^2\big)^2+C\varepsilon_{v}^{1-\frac{1}{2}\rho}\sigma_{v}^{-1}\big(\alpha_{v}^{-2}\sigma_{v}^{-1}A_{v}^2\big)\\
&+2^{v+3}C\varepsilon_{v+1}^{-1}\varepsilon^{(1+\rho)^{v+1}}\\
\leq& C\Big(\varepsilon_{v}^{1-\frac{1}{2}\rho}\big(\alpha_{v}^{-1}\sigma_{v}^{-1}A_{v}\big)^4+\varepsilon_{v+1}\Big)\\
\leq&C\varepsilon_{v}^{1-\frac{35}{36}\rho}\leq\frac{1}{2}
 \end{aligned}
\end{equation}
as $\varepsilon$ small enough. Then $|X_{R^{b_{v+1}}}|^{\lambda_{v}}_{r_{v},D(s_{v}-6\sigma_{v}, r_{v}-4d_{v})\times\mathcal{O}_{v}^{*}}$ is as required.

Moreover, in the same way as (\ref{eqg}), we obtain
\begin{equation}
\begin{aligned}\label{eq3.23}
&|X_{\sum_{n\geq v+2}(\varepsilon^{(1+\rho)^{n}}\tilde{P}^{b_{n}}\circ \Phi_{0} \circ \Phi_{1} \circ \cdots \circ \Phi_{v})}|^{\lambda_{v}}_{r_{v},D(s_{v}-6\sigma_{v},r_{v}-4d_{v})\times\mathcal{O}_{v}^{*}}\\
\leq& 2^{v+3}C\sum_{n\geq v+2}\varepsilon^{(1+\rho)^{n}}|X_{\tilde{P}^{n}}|_{r_{0},D(s_{0},r_{0})\times\mathcal{O}_{v}^{*}}^{\lambda_{0}}\\
\leq&\frac{1}{2}\varepsilon_{v+1}
\end{aligned}
\end{equation}
provided that $\varepsilon$ is small enough. Then (\ref{eqh}) and (\ref{eq3.23}) imply that (\ref{eqi}) is fulfilled with $v+1$ in place of $v$.

We now check that $X_{\hat{P}^{high}_{v+1}}$ satisfy (\ref{Ph}) with $v+1$ in place of $v$.
\begin{equation}\label{eqE}
\begin{aligned}
&|X_{\{\hat{P}_{v}^{high},F_{v}\}^{high}}|^{\lambda_{v}}_{r_{v},D(s_{v}-4\sigma_{v},r_{v}-2d_{v})\times\mathcal{O}_{v}^{*}}\\
\leq& |X_{\{\hat{P}_{v}^{high},F_{v}\}}|^{\lambda_{v}}_{r_{v},D(s_{v}-4\sigma_{v},r_{v}-2d_{v})\times\mathcal{O}_{v}^{*}}
+|X_{\{\hat{P}_{v}^{high},F_{v}\}^{low}}|^{\lambda_{v}}_{r_{v},D(s_{v}-4\sigma_{v},r_{v}-2d_{v})\times\mathcal{O}_{v}^{*}}\\
\leq& C\sigma_{v}^{-1}|X_{\hat{P}_{v}^{high}}|^{\lambda_{v}}_{r_{v},D(s_{v},r_{v})\times\mathcal{O}_{v}^{*}}|X_{F_{v}}|^{\lambda_{v}}_{r_{v},D(s_{v}-3\sigma_{v},r_{v}-d_{v})\times\mathcal{O}_{v}^{*}}\\
\leq& C(\alpha_{v}^{-1}\sigma_{v}^{-1}A_{v})^2|X_{\hat{P}_{v}^{high}}|^{\lambda_{v}}_{r_{v},D(s_{v},r_{v})\times\mathcal{O}_{v}^{*}}.
\end{aligned}
\end{equation}
Recalling (\ref{eqD}), by (\ref{Ph}), (\ref{eqe}), (\ref{eqA}), (\ref{eqB}), (\ref{eqg}), (\ref{eqC}), (\ref{eqE}) and the fact that $$|X_{\tilde{P}^{b_{v+1}}}|^{\lambda_{0}}_{r_{0},D(s_{0},r_{0})\times\mathcal{O}_{v}^{*}}\leq C$$ for $v\geq0$, we obtain

\begin{equation*}
\begin{aligned}
&|X_{\hat{P}_{v+1}^{high}}|^{\lambda_{v}}_{r_{v},D(s_{v}-4\sigma_{v},r_{v}-2d_{v})\times\mathcal{O}_{v}^{*}}\\
\leq&|X_{\hat{P}_{v}^{high}}|^{\lambda_{v}}_{r_{v},D(s_{v},r_{v})\times\mathcal{O}_{v}^{*}}\Big(1+C\varepsilon_{v}\big(\alpha_{v}^{-1}\sigma_{v}^{-1}A_{v}\big)^2\Big)
+C\Big(\varepsilon_{v}^2\big(\alpha_{v}^{-1}\sigma_{v}^{-1}A_{v}\big)^4+\varepsilon_{v+1}\Big)\\
\leq& \big(\varepsilon_{0}+\sum_{j=1}^{v}\varepsilon_{j}^{\frac{1}{2}}\big)(1+\varepsilon_{v+1}^{\frac{1}{2}})+C\varepsilon_{v+1}\\
\leq&\varepsilon_{0}+\sum_{j=1}^{v+1}\varepsilon_{j}^{\frac{1}{2}}.
\end{aligned}
\end{equation*}

 Moreover, (\ref{Bv}) and (\ref{eqj}) imply that
\begin{equation}\label{eq3.30}
\big|\hat{\Omega}_{v}\big|^{\lambda_{v}}_{-\delta,\mathcal{O}_{v+1}^{*}}\leq|X_{B^{b_{v}}}|^{\lambda_{v}}_{r_{v},D(s_{v}-2\sigma_{v},r_{v}-d_{v})\times\mathcal{O}_{v}^{*}}\leq\varepsilon_{v}^{-\frac{1}{8}\rho}.
\end{equation}
Thus, by (\ref{Bv}), the Lipschitz semi-norm of the new frequencies on $\mathcal{O}_{v+1}^{*}$ is bounded by
\begin{equation*}
\begin{aligned}
|\Omega_{v+1}|^{lip}_{-\delta,\mathcal{O}_{v+1}^{*}}&\leq M_{v}+\varepsilon_{v}\big|\hat{\Omega}_{v}\big|^{lip}_{-\delta,\mathcal{O}_{v+1}^{*}}\leq M_{v}+\frac{(M_{v}+1)\varepsilon_{v}}{\alpha_{v}}|X_{B^{b_{v}}}|^{\lambda_{v}}_{r_{v},D(s_{v}-2\sigma_{v},r_{v}-d_{v})\times\mathcal{O}_{v}^{*}}\\
&\leq M_{v}+C(M_{v}+1)\varepsilon_{v}\alpha_{v}^{-2}\sigma_{v}^{-1}A_{v}|X_{R^{b{v}}}|_{r_{v},D(s_{v},r_{v})\times\mathcal{O}_{v}^{*}}^{\lambda_{v}}\\
&\leq (M_{v}+1)(1+\varepsilon_{v}^{1-\frac{1}{4}\rho})-1\\
&\leq (M_{v}+1)(1+2^{-v-1})-1\\
&\leq M_{v+1}
\end{aligned}
\end{equation*}
as required. This completes the proof of the iteration lemma.~~~~$\Box$

\section{Convergence of transformations}

To apply iterative lemma with $v=0$, set $N_{0}=N$, $P_{0}=P$, $s_{0}=s$, $r_{0}=r,\cdots.$ Choosing $\hat{P}_{0}=\varepsilon\tilde{P}^{b_{0}}$, then
(\ref{Ph}) and (\ref{eqi}) with $v=0$ are satisfied by Lemma 2.2.
The small divisor conditions are satisfied by setting $\mathcal{O}_{0}^{*}=\{(\omega^{b_{0}},\omega_{b_{0}}^{'}):\omega^{b_{0}}\in\mathcal{O}_{*}^{b_{0}}\}$, where $\mathcal{O}_{*}^{b_{0}}=[0,1]^{b_{0}}\setminus\bigcup_{k,l}\mathcal{R}_{kl}^{0},$ and $\mathcal{R}^{0}_{kl}=\big\{\omega^{b_{0}}\in[0,1]^{b_{0}}:|\langle k,\omega^{b_{0}}\rangle+\langle l,\Omega_{0}\rangle|< \frac{\alpha_{0}\langle l\rangle_{2}}{(|k|+1)^{2b_{0}+2}}\big\}.$ Hence, using the iterative lemma, we obtain a sequence of transformations $\Phi^{v}= \Phi_{0}\circ\cdots\circ\cdots \Phi_{v-1}:~D_v\times\mathcal{O}_{v-1}^{*}\rightarrow D_{0}$ for $v\geq1,$ and  $H\circ\Phi^{v}=N_{v}+P_{v}$.
We now prove the convergence of $\Phi^{v}$. From (\ref{eq3.33}) and (\ref{eq3.34}), we obtain
\begin{equation}\label{eq3.36}
\begin{aligned}
\frac{1}{\sigma_{v}}|\Phi_{v}-id|^{\lambda_{v}}_{r_{v},D_{v+1}\times\mathcal{O}_{v}^{*}}, |D\Phi_{v}-I|^{\lambda_{v}}_{r_{v},r_{v},D_{v+1}\times\mathcal{O}_{v}^{*}}
\leq C\varepsilon_{v}^{1-\frac{1}{4}\rho}.
\end{aligned}
\end{equation}
We note that the operator norm $|\cdot|_{r,s}$ defined in (\ref{l}) satisfies $|AB|_{r,s}\leq|A|_{r,r}|B|_{s,s}$ for $r\geq s$. For $v\geq1,$ by the chain rule and using (\ref{eq3.36}), we get
\begin{equation*}
\begin{aligned}\
 |D\Phi^{v}|_{r_{0},r_{v},D_{v}\times\mathcal{O}_{v-1}^{*}}\leq\prod_{\mu=0}^{v-1}|D \Phi_{\mu}|_{r_{\mu},r_{\mu},D_{\mu+1}\times\mathcal{O}_{\mu}^{*}}\leq\prod_{\mu=0}^{\infty}(1+\frac{1}{2^{\mu+2}})\leq2
\end{aligned}
\end{equation*}
and

\begin{equation*}
\begin{aligned}\
 |D\Phi^{v}|^{lip}_{r_{0},r_{v},D_{v}\times\mathcal{O}_{v-1}^{*}}&\leq\sum_{\mu=0}^{v-1}|D \Phi_{\mu}|^{lip}_{r_{\mu},r_{\mu},D_{\mu+1}\times\mathcal{O}_{\mu}^{*}}\prod_{0\leq j\leq v-1, j\neq \mu}|D \Phi_{j}|_{r_{j},r_{j},D_{j+1}\times\mathcal{O}_{j}^{*}}\\
 &\leq2\sum_{\mu=0}^{v-1}|D \Phi_{\mu}-I|^{lip}_{r_{\mu},r_{\mu},D_{\mu+1}\times\mathcal{O}_{\mu}^{*}}\leq C\varepsilon_{0}^{1-\frac{11}{36}\rho}
\end{aligned}
\end{equation*}
for sufficiently small $\varepsilon$.

Thus, we have

\begin{equation*}
\begin{aligned}\
|\Phi^{v+1}-\Phi^{v}|_{r_{0},D_{v+1}\times\mathcal{O}_{v}^{*}}
\leq|D\Phi^{v}|_{r_{0},r_{v},D_{v}\times\mathcal{O}_{v-1}^{*}}|\Phi_{v}-id|_{r_{v},D_{v+1}\times\mathcal{O}_{v}^{*}}\leq2|\Phi_{v}-id|_{r_{v},D_{v+1}\times\mathcal{O}_{v}^{*}},
\end{aligned}
\end{equation*}
\begin{equation*}
\begin{aligned}\
|\Phi^{v+1}-\Phi^{v}|^{lip}_{r_{0},D_{v+1}\times\mathcal{O}_{v}^{*}}\leq&|D\Phi^{v}|^{lip}_{r_{0},r_{v},D_{v}\times\mathcal{O}_{v-1}^{*}}|\Phi_{v}-id|_{r_{v},D_{v+1}\times\mathcal{O}_{v}^{*}}\\
&+|D\Phi^{v}|_{r_{0},r_{v},D_{v}\times\mathcal{O}_{v-1}^{*}}|\Phi_{v}-id|^{lip}_{r_{v},D_{v+1}\times\mathcal{O}_{v}^{*}}\\
\leq& C\varepsilon_{0}^{1-\frac{11}{36}\rho}|\Phi_{v}-id|_{r_{v},D_{v+1}\times\mathcal{O}_{v}^{*}}+2|\Phi_{v}-id|^{lip}_{r_{v},D_{v+1}\times\mathcal{O}_{v}^{*}},
\end{aligned}
\end{equation*}
which together with (\ref{eq3.36}) implies that
\begin{equation*}\label{4.46}
\begin{aligned}\
|\Phi^{v+1}-\Phi^{v}|^{\lambda_{0}}_{r_{0},D_{v+1}\times\mathcal{O}_{v}^{*}}\leq\Big(C\varepsilon_{0}^{1-\frac{1}{4}\rho}+2+\frac{2\lambda_{0}}{\lambda_{v}}\Big)|\Phi_{v}-id|^{\lambda_{v}}_{r_{v},D_{v+1}\times\mathcal{O}_{v}^{*}}\leq C\varepsilon_{0}^{\frac{1}{18}\rho}\varepsilon_{v}^{1-\frac{11}{36}\rho}.
\end{aligned}
\end{equation*}
Therefore, the $\Phi^{v}$ converges uniformly on $\bigcap_{v\geq0} \big(D_{v}\times\mathcal{O}_{v}^{*}\big)=D_{*}\times\mathcal{O}^{*}$ to a Lipschitz continuous family of real analytic torus embeddings $\Phi^{\infty}:\mathbb{T}^{\infty}\times\mathcal{O}^{*}\rightarrow \mathcal{P}^{a,p+2},$ where $D_{*}=D(s_{0}/2,r_{0}/2)$ and $\mathcal{O}^{*}=\bigcap_{v\geq0}\mathcal{O}_{v}^{*},$ and

\begin{equation*}\label{4.46}
\begin{aligned}\
|\Phi^{\infty}-id|^{\lambda_{0}}_{r_{0},D_{*}\times\mathcal{O}^{*}}\leq C\varepsilon^{\frac{1}{2}-\frac{1}{8}\rho}.
\end{aligned}
\end{equation*}

\par Thus, at the end of iteration, we obtain the Hamiltonian $H^{\infty}$ of the transformed Hamiltonian system, that is

\begin{equation}\label{}
  H^{\infty}=\langle\omega,J\rangle+\langle \tilde{\Omega},z\bar{z}\rangle+\sum_{|\gamma|_{1}+|\kappa|_{1}\geq3}\hat{P}^{\gamma\kappa}(\theta,\omega)z^{\gamma}\bar{z}^{\kappa},
\end{equation}
where
$\omega=(\omega_{i_{1}},\omega_{i_{2}},\cdots)\in\mathcal{O}^{*},$ $J=(J_{i_{1}},J_{i_{2}},\cdots),$ $i_{j}\in\mathcal{I}_{\infty}$ and $\tilde{\Omega}_{j}$ is close to $\mu_{j}$. It is easy to see that the transformed Hamiltonian system has a solution
\begin{equation*}
  \theta=\omega t+{\rm const.}({\rm mod}~2\pi),~~~~z=\bar{z}=0.
\end{equation*}
Therefore, it is easy to obtain that
for each $\omega=(\omega_{i_{1}},\omega_{i_{2}},\cdots)\in\mathcal{O}^{*},$ the beam equation (\ref{eq1.1})+(\ref{eq1.2}) has an almost-periodic solution of the form
$$u(t,x)=\sum_{j\geq0}\frac{q_{j}(\omega t)\cos(jx)}{\sqrt{\mu_{j}}}$$
where $q_{j}(\omega t),~j=0,1,2\cdots,$ are almost-periodic in $t$ with frequencies $\omega$ and $\|q\|_{a,p+2}=O(\varepsilon^{\frac{1}{2}-\frac{1}{8}\rho})$.\\

\section{ Measure estimate}
At the $v$-th KAM step, we have to exclude the following resonant sets
$$\mathcal{R}^{v}=\bigcup_{k,l}\mathcal{R}_{kl}^{v},$$
where
$$\mathcal{R}_{kl}^{v}=\Big\{\omega^{b_{v}}\in\mathcal{O}^{v}:|\langle k,\omega^{b_{v}}\rangle+\langle l, \Omega_{v}\rangle|<\frac{\alpha_{v}\langle l\rangle_{2}}{(1+v^2)(|k|+1)^{2b_{v}+2}}\Big\},$$

 $$\mathcal{O}^{v}=\{\omega^{b_{v}}:(\omega^{b_{v}},\omega_{b_{v}}^{'})\in\mathcal{O}_{v-1}^{*}\}\subset[0,1]^{b_{v}}$$
with $(k,l)\in\mathcal{Z}^{b_{v}}$ and $\mathcal{O}_{-1}^{*}=\mathcal{O}.$ Here, $\omega^{b_{v}}$ and $\Omega_{v}$ are defined and Lipschitz continuous on $\mathcal{O}_{v-1}^{\ast}$.  Throughout all the iteration steps, we obtain a decreasing sequence of Cantor-like parameter sets $\mathcal{O}\supset\mathcal{O}_{0}^{*}\supset\mathcal{O}_{1}^{*}\supset\cdots$
. Hence, in the limit, we finally get a parameter set $\mathcal{O}^{\ast}=\bigcap_{v=0}^{\infty}\mathcal{O}_{v}^{\ast}.$\\

\indent \textbf {Lemma 5.1}
Let the set $\mathcal{O}=[0,1]^{\infty}$ with probability measure. Then the parameter set $\mathcal{O}^{*}$ obtained above satisfies
\begin{equation*}
 \text{meas}(\mathcal{O}\setminus\mathcal{O}^{*})\leq C \varepsilon^{\frac{1}{48}\rho},
\end{equation*}
where $\text{meas}$ is the standard probability measure on $[0,1]$ and $C>0$ is an absolute constant.\\
\textbf {Proof}
In view of (\ref{eq3.4}) and (\ref{eq3.30}), we can easily get
\begin{equation*}
\begin{aligned}
|\Omega_{v}-\Omega_{0}|_{-\delta}&=\mathop{\text{sup}}\limits_{j\geq0}\Big|\sum_{\tilde{s}=0}^{v-1}\varepsilon_{\tilde{s}}[B_{jj}^{11b_{\tilde{s}}}]\Big|j^{-\delta}
\leq \sum_{\tilde{s}=0}^{v-1}\varepsilon_{\tilde{s}}|X_{B^{b_{\tilde{s}}}}|_{r_{\tilde{s}},D(s_{\tilde{s}},r_{\tilde{s}})\times\mathcal{O}_{\tilde{s}}^{*}}^{\lambda_{\tilde{s}}}\leq C\varepsilon_{0}^{1-\frac{\rho}{8}}<\alpha_{0}=\varepsilon_{0}^{\frac{1}{18}\rho}.
\end{aligned}
\end{equation*}
Moreover, as $\frac{\langle l,\Omega_{0}\rangle}{\langle l\rangle_{2}}\rightarrow 1$ with $\langle l\rangle_{2}\rightarrow\infty$, there exists a positive constant $\beta\gg6\alpha_{0}>0$ such that $|\langle l,\Omega_{0}\rangle|>\beta\langle l\rangle_{2}$.
Thus
$$|\langle l,\Omega_{v}-\Omega_{0}\rangle|\leq|l|_{\delta}|\Omega_{v}-\Omega_{0}|_{-\delta}\leq\langle l\rangle_{2}|\Omega_{v}-\Omega_{0}|_{-\delta}\leq\alpha_{0}\langle l\rangle_{2},$$
where $|l|_{\delta}=\sum|l_{j}|j^{\delta}$ and
$$|\langle l,\Omega_{v}\rangle|>|\langle l,\Omega_{0}\rangle|-|\langle l,\Omega_{v}-\Omega_{0}\rangle|\geq(\beta-\alpha_{0})\langle l\rangle_{2}.$$

\par Case 1. When $|k|\leq\frac{\beta\langle l\rangle_{2}}{4},$
\begin{equation*}
\begin{aligned}
|\langle k, \omega^{b_{v}}\rangle+\langle l,\Omega_{v}\rangle|\geq|\langle l,\Omega_{v}\rangle|-|k||\omega^{b_{v}}|\geq(\beta-\alpha_{0})\langle l\rangle_{2}-\frac{1}{4}\beta\langle l\rangle_{2}>2\alpha_{0}\langle l\rangle_{2}\geq\alpha_{v}\langle l\rangle_{2},
\end{aligned}
\end{equation*}
then $\mathcal{R}_{kl}^{v}$ is empty.
\par Case 2. When $|k|>\frac{\beta\langle l\rangle_{2}}{4}$, let $$g_{0}^{v}(\omega^{b_{v}})=\langle k,\omega^{b_{v}}\rangle,~~~~g_{1}^{v}(\omega^{b_{v}})=\langle k,\omega^{b_{v}}\rangle\pm\Omega_{vj},$$
$$g_{2}^{v}(\omega^{b_{v}})=\langle k,\omega^{b_{v}}\rangle\pm(\Omega_{vi}+\Omega_{vj}),~~~~g_{3}^{v}(\omega^{b_{v}})=\langle k,\omega^{b_{v}}\rangle+\Omega_{vi}-\Omega_{vj}~~(i\neq j),$$
where        $$\Omega_{0j}=\mu_{j}=\sqrt{j^4+m},~~\Omega_{vj}=\sqrt{j^4+m}+O(\varepsilon_{0}^{1-\frac{\rho}{8}}).$$
Choosing a vector $y^{b_{v}}\in\{-1,1\}^{b_{v}}$ such that $\langle k,y^{b_{v}}\rangle=|k|,$ then we obtain
$$\Big|\frac{d}{dt}g_{0}^{v}(\omega^{b_{v}}+ty^{b_{v}})\Big|=|\langle k,y^{b_{v}}\rangle|=|k|>0,$$
$$\Big|\frac{d}{dt}g_{1}^{v}(\omega^{b_{v}}+ty^{b_{v}})\Big|\geq|\langle k,y^{b_{v}}\rangle|-O(\varepsilon_{0}^{1-\frac{\rho}{8}})\geq\frac{1}{3}|k|>0,$$
$$\Big|\frac{d}{dt}g_{2}^{v}(\omega^{b_{v}}+ty^{b_{v}})\Big|\geq|\langle k,y^{b_{v}}\rangle|-O(\varepsilon_{0}^{1-\frac{\rho}{8}})\geq\frac{1}{3}|k|>0,$$
$$\Big|\frac{d}{dt}g_{3}^{v}(\omega^{b_{v}}+ty^{b_{v}})\Big|\geq|\langle k,y^{b_{v}}\rangle|-O(\varepsilon_{0}^{1-\frac{\rho}{8}})\geq\frac{1}{3}|k|>0$$
for sufficiently small $\varepsilon$.
Furthermore,
$$\text{card}\{l~:\langle l\rangle_{2}<\frac{4|k|}{\beta}\}\leq\text{card}\{l~: |l|_{1}<\frac{8|k|}{\beta}\}\leq C \Big(\frac{|k|}{\beta}\Big)^{2}.$$
If we exclude the measure $\mathop{\sum}\limits_{v\geq0}\mathop{\sum}\limits_{0\neq k\in\mathbb{Z}^{b_{v}}}\mathop{\sum}\limits_{\langle l\rangle_{2}<\frac{4|k|}{\beta}}\frac{6\alpha_{v}\langle l\rangle_{2}}{|k|(1+v^2)(|k|+1)^{2b_{v}+2}}$ along some direction, accordingly, exclude the full measure along other directions, then such a residual set is a subset of $\mathcal{O}^{*}$.
First, the excluded measure of the fixed direction satisfies
\begin{equation*}
\begin{aligned}
&\sum_{v\geq0}\sum_{0\neq k\in\mathbb{Z}^{b_{v}}}\sum_{\langle l\rangle_{2}<\frac{4|k|}{\beta}}\frac{6\alpha_{v}\langle l\rangle_{2}}{|k|(1+v^2)(|k|+1)^{2b_{v}+2}}\\
&\leq\sum_{v\geq0}\sum_{0\neq k\in\mathbb{Z}^{b_{v}}}\frac{C\alpha_{v}}{|k|(1+v^2)(|k|+1)^{2b_{v}+2}}\big(\frac{|k|}{\beta}\big)^{3}\\
&\leq\sum_{v\geq0}\frac{C\alpha_{v}}{\beta^3(1+v^2)}\\
&\leq C\alpha_{0}^{\frac{3}{4}}=C \varepsilon^{\frac{1}{48}\rho}.\\
\end{aligned}
\end{equation*}
by the convergence of $\sum_{0\neq k\in\mathbb{Z}^{b_{v}}}\frac{|k|^2}{(|k|+1)^{2b_{v}+2}},$ where $C$ is an absolute constant independent of $v,~\varepsilon.$
Therefore, we have
\begin{equation*}
  \text{meas}(\mathcal{O}\setminus\mathcal{O}^{*})\leq C \varepsilon^{\frac{1}{48}\rho}.
\end{equation*}
The proof of lemma 5.1 is complete.~~~~$\Box$\\

Lemma 5.1 shows that the total measure of all excluded parameter sets can be as small as we wish, and we finally get a Cantor-like parameter set $\mathcal{O}^{*}=\mathop{\bigcap}\limits_{v\geq0}\mathcal{O}_{v}^{*}.$ This completes the proof of Theorem 1.1.\hskip 2.8in $\Box$

\section{Acknowledgements}
This work is supported  by the NNSF(11971163) of China, by Key Laboratory of High Performance Computing and Stochastic Information Processing. The authors would like to thank sincerely Professor Li for friendly suggestions and helpful comments during the preparation of the paper.\\
\\

\section{Appendix}
\indent \textbf {Lemma 7.1} [3]
 For $\sigma>0$ and $v>0,$ the following inequalities hold true:
\begin{equation*}
\begin{aligned}\
\sum_{k\in\mathbb{Z}^{n}}e^{-2|k|\sigma}\leq\frac{1}{\sigma^{n}}(1+e)^{n},
\end{aligned}
\end{equation*}
\begin{equation*}
\begin{aligned}\
\sum_{k\in\mathbb{Z}^{n}}e^{-2|k|\sigma}|k|^{v}\leq\Big(\frac{v}{e}\Big)^{v}\frac{1}{\sigma^{v+n}}(1+e)^{n}.
\end{aligned}
\end{equation*}
\\
\indent \textbf {Lemma 7.2} [31]
If $A=(A_{ij})$ is a bounded linear operator on $\ell^2,$ then also $B=(B_{ij})$ with
$$B_{ij}=\frac{|A_{ij}|}{|i-j|},~~~~~~i\neq j,$$
and $B_{ii}=0$ is a bounded linear operator on $\ell^2,$ and $\|B\|\leq\frac{\pi}{\sqrt{3}}\|A\|$.
\\
\indent \textbf {Lemma 7.3}\label{lem5.3}
For $k\in\mathbb{Z}^{b_{v}}$, we have that
$$\sqrt{2^{b_{v}}\sum_{k\in\mathbb{Z}^{b_{v}}}|k|^2[(1+v^2)(|k|+1)^{2b_{v}+2}]^4e^{-2|k|\sigma_{v}}}\leq\Big(\frac{16(2b_{v}+3)}{e}\Big)^{4b_{v}+6}\frac{1}{\sigma_{v}^{5b_{v}+6}}.$$

\textbf{Proof}
 Since $2^{9b_{v}+8}(1+v^2)^4\cdot4^{b_{v}}\leq4^{8b_{v}+12},$ and using Lemma 7.1, we obtain
\begin{equation*}
 \begin{aligned}
 &2^{b_{v}}\sum_{k\in\mathbb{Z}^{b_{v}}}|k|^4[(1+v^2)(|k|+1)^{2b_{v}+2}]^4e^{-2|k|\sigma_{v}}\\
 &\leq2^{b_{v}}(1+v^2)^4\sum_{k\in\mathbb{Z}^{b_{v}}}|k|^4(|k|+1)^{8b_{v}+8}e^{-2|k|\sigma_{v}}\\
 &\leq2^{9b_{v}+8}(1+v^2)^4\Big(\frac{8b_{v}+12}{e}\Big)^{8b_{v}+12}\frac{1}{\sigma_{v}^{9b_{v}+12}}(1+e)^{b_{v}}\\
&\leq2^{9b_{v}+8}(1+v^2)^4\Big(\frac{8b_{v}+12}{e}\Big)^{8b_{v}+12}\frac{1}{\sigma_{v}^{9b_{v}+12}}4^{b_{v}}\\
&\leq\Big(\frac{16(2b_{v}+3)}{e}\Big)^{8b_{v}+12}\frac{1}{\sigma_{v}^{9b_{v}+12}}.\\
\end{aligned}
 \end{equation*}
Thus
$$\sqrt{2^{b_{v}}\sum_{k\in\mathbb{Z}^{b_{v}}}|k|^2[(1+v^2)(|k|+1)^{2b_{v}+2}]^4e^{-2|k|\sigma_{v}}}\leq\Big(\frac{16(2b_{v}+3)}{e}\Big)^{4b_{v}+6}\frac{1}{\sigma_{v}^{5b_{v}+6}}.~~\Box$$
\\
\textbf{Data Availability}

The data that supports the findings of this study are available within this article.

{}

\end{document}